\documentclass[leqno,oneside,english]{article}
\usepackage[latin1]{inputenc}
\usepackage{amsmath}
\usepackage{amssymb}
\usepackage{amsthm}
\usepackage{color}
\usepackage{dsfont}
\usepackage{tikz}
\usepackage{hyperref}
\usepackage{enumerate}
\usepackage{eucal}
\numberwithin{equation}{section}
\usepackage{subfigure}
\usepackage{pgfplots}

\newcommand{\madeline}[1]{\textcolor{black}{#1}}

\makeatletter
\usepackage[T1]{fontenc}
\usepackage[latin1]{inputenc}

\usepackage{graphicx}
\def\qed{\hbox{\hskip 6pt\vrule width6pt 
height7pt
depth1pt  \hskip1pt}\bigskip}

\let\oldqedhere\qedhere
\renewcommand{\qedhere}{\pushQED{\qed}\oldqedhere}

\newcommand{\re}{\mathrm{Re}}

\newcommand{\im}{\mathrm{Im}}

\definecolor{monvert}{rgb}{0,0.5,0}
\definecolor{monrouge}{rgb}{0.6,0,0}
\definecolor{monbrun}{rgb}{0.5,0,0.5}
\definecolor{monbleu}{rgb}{0,0,0.5}
\makeatletter \textwidth15cm \textheight20truecm

\newcommand{\R}{{\mathbb R}}
\newcommand{\C}{{\mathbb C}}
\newcommand{\N}{{\mathbb N}}

\def\be{\begin{equation}}
\def\ee{\end{equation}}
\def\beq{\begin{eqnarray}}
\def\eeq{\end{eqnarray}}
\def\beqs{\begin{eqnarray*}}
\def\eeqs{\end{eqnarray*}}
\def\ea{\end{array}}
\def\ea{\end{array}}
\def\bt{\begin{thm}}
\def\et{\end{thm}}
\def\br{\begin{rk}}
\def\er{\end{rk}}
\def\bl{\begin{lem}}
\def\el{\end{lem}}
\def\bc{\begin{cor}}
\def\ec{\end{cor}}
\def\bex{\begin{exo}}
\def\eex{\end{exo}}
\renewcommand{\div}{\mathop{\rm div}\nolimits}

\newtheorem{thm}{Theorem}[section]

\newtheorem{cor}[thm]{Corollary}

\newtheorem{defi}[thm]{Definition}
\newtheorem{exo}[thm]{Example}

\newtheorem{lem}[thm]{Lemma}

\newtheorem{rk}[thm]{Remark}

\newcommand{\Om}{\Omega}

\makeatother

\usepackage{babel}
\author{Rolando Magnanini
\footnote{Universit\`a degli Studi di Firenze, Dipartimento di Matematica e Informatica U. Dini,
Viale Morgagni 67/A, 
50134 Firenze, {Italy,}
rolando.magnanini@unifi.it}, Serge Nicaise
\footnote{Universit\'e Polytechnique Hauts-de-France, C\'ERAMATHS/DMATHS and FR CNRS 2037,
F-59313 - Valenciennes Cedex 9 France,
Serge.Nicaise@uphf.fr},   Madeline Chauvier \footnote{Universit\'e Polytechnique Hauts-de-France, C\'ERAMATHS/DMATHS and FR CNRS 2037,
F-59313 - Valenciennes Cedex 9 France and
  D\'epartement de Math\'ematique,
  Universit\'e de Mons,
  place du parc~20,
  B-7000 Mons, Belgium, Madeline.CHAUVIER@umons.ac.be}}
\begin{document}
\title{Critical points of solutions of elliptic equations in divergence form  in planar non simply connected domains with smooth or  nonsmooth   boundary}

  \maketitle

\begin{abstract}
 We study   the critical points of  the solution of  
second elliptic equations      in divergence and diagonal form
with a bounded and positive \madeline{definite} coefficient,
under the assumption that the statement of the Hopf lemma
holds  (sign assumptions on its normal derivatives) along the boundary.
 The proof 
combines the argument principle introduced in \cite{AlessandriniMagnanini:92} for elliptic equations
 with the representation formula (using quasi-conformal mappings) for operators in divergence form in simply connected domains 
     \cite{AlessandriniMagnanini:94}.  The case of a degenerate coefficient is also treated where we combine the  level lines technique and the maximum principle
with the argument principle.  Finally, some numerical experiments on illustrative examples are presented. 
\end{abstract}

\noindent{\bf  AMS (MOS) subject classification: 35J25, 35B38, 58E05}  

\noindent{\bf Key Words:} Critical points, level lines, elliptic equations in divergence form
\vspace{1cm}

\section{Introduction}

\begin{figure}
\begin{center}
\begin{tikzpicture}[scale = 0.9,baseline={([yshift=-.6ex]current bounding box.center)}]
            \draw (0:0) arc (0:180:3cm);
            \draw (-6,0) -- (0,0);
            \draw (-3,1.3) circle (0.5);
            \draw (-3,0.5) node {$\Gamma_c$};
            \draw (-2, 3.1) node {$\Gamma_d$};
            \draw (-0.5,0.22) node {$\Gamma_d$};
        \end{tikzpicture}
        \end{center}
\caption{An illustrative configuration\label{exempleconfig}} 
\end{figure}
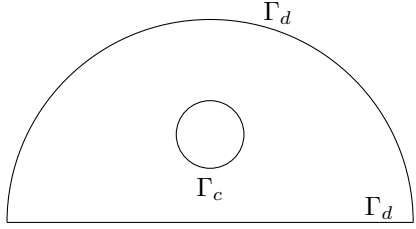

 Recently    the engineering community  has  
an  increasing  interest  of the transmission of electricity using direct current, because the production of this type of current  grows due to the increase of the production of the renewable energies and  the direct current is more suitable for the long distance transport.  The simplest mathematical model used by this community \cite{SharmaJaniszewski:69,ButlerCendesHoburg:89,CristinaDinelliFeliziani:91,Lobry:14,AmorusoLattarulo:14} is
the ion flow field in the ambient of conductors, that  is actually governed by Poisson's equation and
the current density conservation equation set in a region $\Omega$ of the plane:
 \begin{equation}\label{sys : gen}
            \begin{cases}
          -\Delta \varphi= \rho , & \hbox{ in } \Omega,\vspace{2mm}  \\
            \div \bigl(\rho  \nabla\varphi \bigr) = 0, & \hbox{ in } \Omega, 
            \end{cases} 
        \end{equation}   
where the variable $\varphi$ represents the electric potential and the variable $\rho$ the space charge density (that for physical reasons  has to be non negative) are the two unknowns.
The domain    $\Omega$ corresponds to the air  surrounding the conductors  
above the ground  and should be an unbounded domain, but
  for computational reasons it is reduced to a bounded one.
  It is admitted \cite{SharmaJaniszewski:69,ButlerCendesHoburg:89,CristinaDinelliFeliziani:91,Lobry:14,AmorusoLattarulo:14} that this   system of nonlinear partial differential equations is completed with the  {Dirichlet} boundary conditions
   \begin{equation} \label{sys : BC}
                 \begin{cases}
                    \varphi = 1, &  \hbox{ on }\Gamma_c,\\
                    \varphi = 0,& \hbox{ on } \Gamma_d, 
                \end{cases}
            \end{equation}
{and the Neumann condition
\begin{equation} \label{neumann}
\frac{\partial \varphi}{\partial n} = E_0, \hbox{ on } \Gamma_c,
\end{equation}
where $\frac{\partial\varphi}{\partial n}$ denotes the exterior normal derivative} and $E_0$ is a positive constant.  Here, $\Gamma_c$ represents the boundary of the  conductor and $\Gamma_d$ represents the truncated boundary {between air and ground. The two boundaries} are such that $\partial \Omega = \Gamma_c \cup \Gamma_d$ and $\Gamma_c \cap \Gamma_d = \emptyset$, see Figure \ref{exempleconfig} for an illustrative configuration.

This problem    
is studied in  \cite{Chauvieretco:2025}, where in its full generality  the  existence of a solution {of \eqref{sys : gen},\eqref{sys : BC}}
is proved, some properties of the solution are given, and a numerical algorithm based on the minimisation of a functional  is proposed to find an approximation of such a solution. 
In the engineering community,  another algorithm is proposed to find an approximated solution. It is based on the resolution  {of the first equation} from \eqref{sys : gen} by a finite element method
and then on the resolution of the  {second equation by} using the method of characteristic (see \cite{SharmaJaniszewski:69,ButlerCendesHoburg:89,GuillodPfeifferFranck:14,Xiaoetall:17,Zhangetall:15}), transforming it as follows:
\be\label{eq:transport}
0= \div \bigl(\rho  \nabla\varphi \bigr) = \nabla \varphi\cdot\nabla \rho+ \rho \Delta 
 \varphi=\nabla \varphi\cdot\nabla \rho-\rho^2.
 \ee
But such a method can  be implemented only  if  the set of points $v\in \bar\Omega$ where  $\nabla \varphi(v)$ is  zero is finite.
The goal of the present  paper is therefore to give the exact value of
the number of critical points of $\varphi$ (i.e. points $v$ such that $\nabla \varphi(v)=(0,0)^\top$)   for some domains $\Omega$ and under  appropriate and natural 
 conditions on $\rho$, namely either $\rho\in L^\infty(\Omega)$ is uniformly positive or
 is uniformly positive except near the  corners of the exterior boundary of the domain where it can tend to zero.
 This 
degenerate case is completely natural    if the exterior boundary has a corner $v$,
because by the second boundary condition from \eqref{sys : BC}, 
$\nabla \varphi(v)=(0,0)^\top$, assuming that $\varphi$ is smooth enough.
Hence  the transport equation \eqref{eq:transport}  implies that
$\rho$ tends to zero along a characteristic line that hits a corner.

Our main results are based on the study of the critical points of  the solution of the 
second equation from \eqref{sys : gen},
\begin{equation}\label{eq:divrhonabla}
\div \bigl(\rho  \nabla\varphi \bigr) = 0,   \hbox{ in } \Omega, 
\end{equation}
subject to the Dirichlet boundary conditions
   \begin{equation} \label{sys : DirBC}
                \begin{cases}
                    \varphi = 1, &  \hbox{ on }\Gamma_c,\\
                    \varphi = 0,& \hbox{ on } \Gamma_d,
                \end{cases}
            \end{equation}
and under the assumption that the statement of the Hopf lemma
holds  (sign assumptions on its normal derivatives) along $\partial \Omega$. 
 Their proofs actually
combine the argument principle introduced in \cite{AlessandriniMagnanini:92} for {elliptic equations} 
 with the representation formula (using quasi-conformal mappings) for operators in divergence form in simply connected domains
     \cite{AlessandriniMagnanini:94}, allowing to characterize the level lines of the solutions of 
\eqref{eq:divrhonabla}. Such   results are already known 
for the Laplace equation (corresponding to $\rho=1$) \cite{AlessandriniMagnanini:92},
for linear or quasilinear second order operators in non-divergence form  \cite{Alessandrini:87, AlessandriniLupoRosset:93, Dengetco:18,Dengetco:22}, for  $\rho \in W^{1,\infty}(\Omega)$ \cite{Alessandrini:87b}, or 
$\rho\in L^\infty(\Omega)$ but for $\Omega$ simply connected  \cite{AlessandriniMagnanini:94}.  However, to the best of our knowledge, the case of non simply connected domains and 
 $\rho$ only in $L^\infty(\Omega)$ (and even degenerate) is not yet treated.

Let us further notice that combining some conformal mapping transformation and the level lines technique,  we also prove that the Hopf lemma
holds for  the solution of \eqref{eq:divrhonabla}  with  the Dirichlet boundary conditions
     \eqref{sys : DirBC} under some regularity assumptions on the  boundary of $\Omega$
     and the assumption that $\rho\in L^\infty(\Omega)$ is uniformly bounded from below.
     Such a result  holds if $\partial \Omega$ is smooth enough and $\rho\in C^{0,\alpha}(\bar \Omega)$, with $\alpha\in (0,1]$, see \cite[Lemma 7]{FinnGilbarg:57} or \cite[Theorem 1.1]{SabinadeLis:15}. In the opposite case,
   it can be wrong, but the counterexamples from the literature recalled 
   in  \cite{SabinadeLis:15} do not enter in our setting.

The paper is organized as follows: in Section \ref{s:setting}
we present the problem and formulate the required  assumptions on its solution, namely that   Hopf lemma holds on  the exterior/interior boundary.
 In Section \ref{s:criticalpts}, using local quasi-conformal changes of variables, we show that the critical points of the solution are isolated.
Based on the argument principle, we  prove in  Section  \ref{s:numbercriticalpts} that  the number of  critical points inside the domain is equal to the number of connected components of the obstacles minus 1.
Section \ref{s:degeneratedcase} is devoted to the extension of our results to the case of 
degenerate density, there we combine the  level lines technique and the maximum principle
with the argument principle.
In Section \ref{s:Hopf}, under some regularity assumptions on the exterior/interior boundary,
we prove   the Hopf lemma on the exterior/interior boundary with the help  of level lines 
technique. Finally in  Section \ref{sexample} some numerical tests on illustrative examples are presented. 
 Since the main concerns of this paper are the characterization of the critical points, the proof  that conformal mappings preserve  elliptic equations in divergence and diagonal form
  is postponed to Appendix \ref{sappendixA}.

Let us finish this introduction with some notation used in the
remainder of the paper:   The usual
Sobolev space of order $s>0$ in an open set $\Omega$ of the plane is denoted by $H^{s}(\Omega)$.
The set $\N$ (resp. $\N^*$) is the set of non negative (resp. positive) integers.

 \section{Setting of the problem and assumptions\label{s:setting}}
\hspace{5 mm}
\setcounter{equation}{0}

Let $\mathcal{O}\subset\R^2$  be a bounded connected open set and {$\mathcal{B}$ \madeline{be} the finite union of open simply connected domains $B_1, \dots, B_M\subset\R^2$, with pairwise disjoint closures.} 
Suppose that $\bar{\mathcal{B}} \subset \mathcal{O}$ and set $\Omega=\mathcal{O}\setminus \bar{\mathcal{B}}$. It is clear that $\Omega$ is connected. We assume that 
the boundary of $\mathcal{B}$ is $C^{1,1}$. \madeline{We further assume
that the boundary of $\mathcal{O}$ is smooth except at a finite number of points,  {which we call vertices,} in the following sense:}

\medskip\noindent
{\bf (A1)} either the boundary of $\mathcal{O}$ is $C^{1,1}$,  \madeline{in which case we set $I=0$,} or the boundary of $\mathcal{O}$ is $C^{1,1}$, except at a finite number of points
$v_i\in \partial \mathcal{O}$, $i=1,\cdots, I$, with $I\in \N^*$. \madeline{In that last case}, the boundary of $\Omega$ coincides locally with the image by a conformal map of  a quarter plane (\madeline{when we identify} $\R^2$ with $\C$). Namely there exist
a neighbourhood $U_i$ of $v_i$ in $\C$,
an open, non empty, simply connected and bounded set $V_i$ of $\C$ such that $0\in V_i$
and a conformal 
mapping 
\[
\Phi_i:V_i\cap Q\to U_i\cap \Omega,
\]
 where
\[
{Q:=\{w\in \C: \re(w)>0,\  \im(w)>0\}.}
\]
By convention, if the boundary of $\mathcal{O}$ is $C^{1,1}$, we take $I=0$.

We \madeline{assume that the conductivity} $\rho\in L^\infty(\Omega)$ \madeline{satisfies the bound}
\be\label{bdrho}
\delta \leq \rho \hbox{ a.e. in } \Omega,
\ee
\madeline{for some positive real number} $\delta$.

Now we consider the unique weak solution $u\in H^1(\Omega)$ of
\be\label{EDPinu}
\div(\rho \nabla u)= 0\hbox{ in } \Omega,
\ee
with the boundary conditions
\be\label{bcu}
\begin{cases}
u= 1\hbox{ on } \partial \mathcal{B},
\\
u= 0\hbox{ on } \partial \mathcal{O}.
\end{cases}
\ee

In the whole paper, we assume that this solution 
belongs to $C^1(\bar \Omega)$. This assumption is motivated by the following
regularity property of a solution $\varphi\in H^1(\Omega)$ of   \eqref{sys : gen}
with the   Dirichlet boundary condition \eqref{sys : DirBC}. Indeed 
 thanks to elliptic regularity results,   $\varphi\in H^1(\Omega)$ \madeline{is} also a solution of  the first equation in \eqref{sys : gen}, with the   Dirichlet boundary condition \eqref{sys : DirBC}. \madeline{Thus, $\varphi$}  turns out to be in 
$W^{2,p}(\Omega)$, for all $p\geq 2$, due to the assumptions on $\Omega$ and $\rho$. Hence, $\varphi\in C^1(\bar \Omega)$.

{
\br\label{rk:reg}
{\rm Using a smooth lifting for the boundary condition on \madeline{$\partial \mathcal{B}$}, by \cite[Theorem 1] {Meyers:63}, the solution $u\in H^1(\Omega)$ of  
\eqref{EDPinu}-\eqref{bcu}  belongs to 
$W^{1,p}(\Omega)$ for some $p>2$, hence  belongs to $C^{0, \alpha}(\bar \Omega)$, for some $\alpha\in (0,1)$ (see  \cite{PiccininiSpagnolo:72} for an optimal interior regularity result). If $\rho\in W^{1,q}(\Omega)$ for some $q\geq 2$ and $\Omega$ has a $C^2$ boundary, some $W^{2,p}(\Omega)$ estimates, with $p>1$
are also available in \cite{MR3500302,Perelmuter:25}. However, when the coefficient $\rho$ shows some simple jump discontinuities, the solution is expected to be Lipschitz continuous. 
Indeed if the discontinuities are made of a finite union of disjoint smooth curves that do not intersect the   boundary and the coefficient $\rho$ is smooth on each relevant subdomain, then this solution is piecewise in $W^{2,p}$ with $p>2$ (see \cite{Seftel:63} for instance). \madeline{Thus, the solution} is piecewise $C^1$ (in the sense that the gradient is continuous except across the discontinuity curves), \madeline{and in particular}  is Lipschitz continuous. 
\madeline{For instance,} in \cite{Magnanini:80}, the first author of this paper explicitly computed the series solution of the Dirichlet problem in a ball in the case in which the coefficient is radially symmetric and has a finite number of jump discontinuities
 and showed that the solution is indeed Lipschitz continuous. For piecewise smooth coefficients that may cross or may intersect the  boundary, some regularity results in the form of a decomposition into a regular part that is piecewise in $W^{2,p}$, with $p>2$, and explicit singular functions are also available in \cite{nicaise:94a}. In such a case, the solution is still piecewise $C^1$ except eventually at the crossing points, \madeline{but} it is no more Lipschitz continuous in general.}
 \er
 }

As $\Omega$ is connected, due to the boundary conditions \eqref{bcu}, $u$ cannot be constant and therefore  the  strong maximum principle (see \cite[Theorems 8.1 and 8.19]{gilbarg:77}
for instance) yields that
\be\label{bounds on u}
0<u<1 \hbox{ in } \Omega.
\ee
This property implies that 
\beqs
\partial_n u\leq 0 \hbox{ on } \partial \mathcal{O},
\\
\partial_n u\geq 0 \hbox{ on } \partial \mathcal{B},
\eeqs
but we will assume more. Namely we assume that the statement of Hopf lemma
holds on the sets $\partial \mathcal{O}\setminus \cup_{i=1}^I \{v_i\}$
\footnote{If  $\partial \mathcal{O}$ is $C^{1,1}$, then $\mathcal{O}$ has no corner, and therefore $\cup_{i=1}^0 \{v_i\}=\emptyset$} and $\partial \madeline{\mathcal{B}}$:
 \beq
 \label{assumderiveenormaleextbdy}
&&\partial_n u< 0 \hbox{ on } \partial \mathcal{O}\setminus \cup_{i=1}^I \{v_i\},
\\
&&\partial_n u>0 \hbox{ on } \partial \mathcal{B}.
 \label{assumderiveenormaleintbdy}
\eeq

\br{\rm
For the applications to system  \eqref{sys : gen}, the  assumptions  \eqref{assumderiveenormaleextbdy} and \eqref{assumderiveenormaleintbdy} 
are natural.
Indeed for a  solution $\varphi\in C^1(\bar \Omega)$ of \eqref{sys : gen} and \eqref{sys : DirBC}, \madeline{condition} \eqref{assumderiveenormaleextbdy} holds by applying \cite[Lemma 3.4]{gilbarg:77}
to the first equation of  \eqref{sys : gen} with the boundary conditions \eqref{sys : DirBC}. On the contrary \eqref{assumderiveenormaleintbdy} holds by applying 
\cite[Lemma 7]{FinnGilbarg:57} or \cite[Theorem 1.1]{SabinadeLis:15}
to the second equation of  \eqref{sys : gen} with the boundary conditions \eqref{sys : DirBC} under the assumption that 
$\rho\in C^{0,\alpha}(\bar \Omega)$, with $\alpha\in (0,1]$. We shall show in Section \ref{s:Hopf}
that these properties remain valid under some regularity assumptions on $\partial \mathcal{O}$
and $\partial \mathcal{B}$.
}
\er

\br{\rm
\madeline{The results below remain valid if the boundary conditions \eqref{bcu} are replaced by 
\[
\begin{cases}
u= a\hbox{ on } \partial \mathcal{B},
\\
u= b\hbox{ on } \partial \mathcal{O},
\end{cases}
\]
where $a,b$ are two positive real numbers such that $a>b$. Indeed, it suffices to apply our results to $\tfrac{u-b}{a-b}$ that satisfies \eqref{EDPinu} and \eqref{bcu}.}
}\er

\section{The critical points are isolated \label{s:criticalpts}}

Let us recall the next definition.
\begin{defi}
A point $x\in \bar \Omega$ is a (geometric) critical point of 
$u\in C^1(\bar \Omega)$ if and only if
\[
\nabla u(x)=0.
\]
\end{defi}

\subsection{Case of a $C^{1,1}$ exterior boundary}

\bl\label{l:ptcritiqueregularcase}
Assume that the solution $u$ of
\eqref{EDPinu}-\eqref{bcu} 
belongs to $C^1(\bar \Omega)$
and satisfies \eqref{assumderiveenormaleextbdy} and 
 \eqref{assumderiveenormaleintbdy}. \madeline{If} $\partial \mathcal{O}$ is $C^{1,1}$, \madeline{then} $u$ has   a finite number
of  critical points in $\Omega$ (and no critical point on $\partial \Omega$).
\el
\begin{proof}
Due to the assumptions \eqref{assumderiveenormaleextbdy} and
 \eqref{assumderiveenormaleintbdy} and the regularity $u\in C^1(\bar \Omega)$, there exists
 a neighbourhood 
$W$ of $\partial \Omega$ such that
\be\label{gradientnonzero}
\nabla u\ne (0,0)^\top \hbox{ on } W\cap \bar \Omega.
\ee

Now we cover $\Omega$ by overlapping patches $\Omega_j$, $j=1,\cdots, J$
that are  simply connected with a Lipschitz boundary, 
and such that  {the surface measures of 
$\Omega_j\cap \partial \mathcal{B}$ and $\Omega_j\cap \partial \mathcal{O}$ are positive},
for all $j=1,\cdots, J$. Then
we can apply \madeline{\cite[Lemma 2.5 and the following remark]{AlessandriniMagnanini:94}} to $u$ restricted to $\Omega_j$ (since $u$ cannot be constant) to deduce that the critical points of  
$u$ in $\Omega_j$ are isolated without accumulation point inside $\Omega_j$
(they could accumulate near the boundary of $\Omega_j$). 

We conclude by  observing that the $\Omega_j$'s overlap among them and with
$W\cap \bar \Omega$.
\end{proof}

\subsection{Case of a piecewise $C^{1,1}$ exterior boundary}

\bl\label{l:ptcritiquenonregularcase}
Assume that the solution $u$ of
\eqref{EDPinu}-\eqref{bcu} 
belongs to $C^1(\bar \Omega)$
and satisfies \eqref{assumderiveenormaleextbdy} and 
 \eqref{assumderiveenormaleintbdy}. \madeline{If} {\bf (A1)} holds with $I\geq 1$, \madeline{then} $u$ has a finite number
of  critical points in $\Omega$, and the sole   critical points of $u$ on $\partial \Omega$
are the vertices $v_i$, $i=1,\cdots, I$.
\el
\begin{proof}
Under the assumption {\bf (A1)} with $I\geq 1$, the arguments of the proof of Lemma \ref{l:ptcritiqueregularcase} only fail in a neighbourhood of $v_i$, $i=1,\cdots, I$.
Notice that {\bf (A1)} with $I\geq 1$ and the second boundary condition from \eqref{bcu} yield that
\[
\nabla u(v_i)=(0,0)^\top, \forall i=1,\cdots, I.
\]

Hence for any $i=1,\cdots, I$,  we only have to show that there exists  $\delta_i>0$
such that $u$ has only $v_i$ as critical point in $\bar\Omega\cap B(\madeline{v_i},\delta_i)$.
For that purpose, we fix $i=1,\cdots, I$ and drop the index $i$ for shortness.
Then we proceed in three steps:
\\
{\bf Step 1: Straighten the edges:}
Using the mapping $\Phi$ from the assumption  {\bf (A1)}, we denote by
\[
\rho^*(w)=\rho(\Phi(w)),  \forall w\in V^*:=V\cap   Q,
\]
and
\[
u^*(w)=u(\Phi(w)), \forall w\in  \overline{V^*}.
\]
By Lemma \ref{l:preservation}, $u^*\in C^1(\overline{V^*})$ is \madeline{a} solution of
\be\label{pb*}
\begin{cases}
\div(\rho^* \nabla u^*)= 0\hbox{ in } V^*,
\\
u^*= 0\hbox{ on } \partial V^*\cap \partial Q.
\end{cases}
\ee
\\
{\bf Step 2: Odd reflection with respect to the axis $\madeline{\re(w)}=0$:}
With the notation $w\equiv (x_1,y_1)$, with $x_1,y_1\in \R$, we set
\[
V^{**}=\{(x_1,y_1)\in \R^2\,:\, (|x_1|,y_1)\in V^*\}
\]
and 
\beqs
\rho^{**}(x_1,y_1)=
\begin{cases}
 \rho^*(x_1,y_1) &\hbox{ if } (x_1,y_1)\in V^*,
\\
\rho^*(-x_1,y_1) &\hbox{ else, }
\end{cases}
\\
u^{**}(x_1,y_1)=
\begin{cases}
 u^*(x_1,y_1) &\hbox{ if } (x_1,y_1)\in \overline{V^*},
\\
-u^*(-x_1,y_1) &\hbox{ if } (x_1,y_1)\in \overline{V^{**}}\setminus \overline{V^*}.
\end{cases}
\eeqs
We directly check that $u^{**}\in H^1(V^{**})\cap C(\overline{V^{**}})$ is \madeline{a} solution of
\be\label{pb**}
\begin{cases}
\div(\rho^{**} \nabla u^{**})= 0 &\hbox{ in } V^{**},
\\
u^{**}= 0 &\hbox{ on } \partial V^{**}\cap \partial \pi_+,
\end{cases}
\ee
where $\pi_+:=\{(x_1,y_1)\in \R^2: \madeline{y_1}>0\}$ is the  right-hand  half-plane.
\\
 {\bf Step 3: Odd reflection with respect to the axis $\madeline{\im(w)}=0$:}
As before, by setting
\[
V^{***}=\{(x_1,y_1)\in \R^2\,:\, (|x_1|,|y_1|)\in V^*\},
\]
and 
\beqs
\rho^{***}(x_1,y_1)=
\begin{cases}
 \rho^{**}(x_1,y_1) &\hbox{ if } (x_1,y_1)\in V^{**},
\\
\rho^{**}(x_1,-y_1) &\hbox{ else, }
\end{cases}
\\
u^{***}(x_1,y_1)=
\begin{cases}
 u^{**}(x_1,y_1) &\hbox{ if } (x_1,y_1)\in \overline{V^{**}},
\\
-u^{**}(x_1,-y_1) &\hbox{ if } (x_1,y_1)\in \overline{V^{***}}\setminus \overline{V^{**}},
\end{cases}
\eeqs
we find that $u^{***}\in H^1(V^{***}) \madeline{\cap  C(\overline{V^{***}})}$  and is a solution of
\be\label{EDP***}
\div(\rho^{***} \nabla u^{***})= 0\hbox{ in } V^{***}.
\ee
 
 Note that by construction 	and the second boundary condition from \eqref{bcu},
 $u^{***}=0$ on the segments $S_1:= V^{***}\cap \{ x_1=0\}$ and $S_2=V^{***}\cap \{ y_1=0\}$, while due to \eqref{bounds on u},
 \be\label{eq:12}
0<|u^{***}|<1 \hbox{ on } \overline{ V^{***}}\setminus (S_1\cup S_2).
\ee
Applying \cite[Theorem 2.1]{AlessandriniMagnanini:94} to $u^{***}$ in $V^{***}$, we deduce  that the critical points of  
$u^{***}$ in $V^{***}$ are isolated without accumulation points inside $V^{***}$.
Hence in a sufficiently small neighbourhood $W^{***}$ of $0$ (in $\R^2$), the sole critical point of 
$u^{***}$ in $V^{***}\cap W^{***}$ is $0$.

Coming back to $u$,  there are no other critical points of $u$ near $v$.
\end{proof}

\section{Number of critical points inside the domain\label{s:numbercriticalpts}}

The goal of this section is to prove that, under our original assumptions, the 
number  \madeline{of  critical points} in $\Omega$ of the solution $u$ of problem  
\eqref{EDPinu}-\eqref{bcu}  is exactly the number of connected components of $\madeline{\mathcal{B}}$ minus 1. Our proof  is based on the argument principle. 


Let us start with the notion of index.
\begin{defi}
Let $u\in C^1(\bar \Omega)$ be a solution of
\eqref{EDPinu}-\eqref{bcu}.
Let $\Omega'\subset \Omega$ be a simply connected domain with a Lipschitz boundary, and
let $D\subset \subset \Omega'$. If $u$ has no geometric critical points 
on $\partial D$, then  
the index of $u$ in $D$ is defined by
\be\label{defIu}
I(D,u)=I(\chi(D), h)=-\frac{1}{2\pi}\int_{\madeline{+}\partial \chi(D)} d\arg (\nabla h),
\ee
where $\chi$ is the quasi-conformal mapping from $\Omega'$ to the unit ball $B_1(0)$
built up in \cite[Theorem 2.1]{AlessandriniMagnanini:94} such that
$h=u\circ \chi^{-1}$ is harmonic in $B_1(0)$. \madeline{In other words, $-2\pi I(D,u)$ is the increase of the angle 
$\arg(\nabla h)$  along the contour $+\partial\chi(D)$ oriented counterclockwise.}

Moreover, we define the geometric index of $u$ at $z_0\in \Omega$ as
\be\label{defIuatz0}
I(z_0,u)=\lim_{r\to 0} I(B_r(z_0), u).
\ee
\end{defi}

As $h$ is harmonic, $I(B_1(0), h)$ is the number of critical points of $h$ inside $B_1(0)$
(when counted according to their multiplicities), and hence the number  critical points of $u$ inside $D$.

It turns out that   the $C^1$ regularity of $u$ leads to the invariance of the index, namely
the next result holds.
\bl\label{l:indexu}
Let $u\in C^1(\bar \Omega)$ be a solution of
\eqref{EDPinu}-\eqref{bcu}.
Let $\Omega'\subset  \Omega$ be a simply connected domain with a Lipschitz boundary, and
let $D\subset \subset \Omega'$. Assume further that $u$ has no geometric critical points 
on $\partial D$, then 
\be\label{defIulocal}
I(D,u)=-\frac{1}{2\pi}\int_{+\partial D} d\arg (\nabla u).
\ee
\el
\begin{proof}
  \madeline{
  As quasi-conformal mappings are
  orientation-preserving mappings, the quasi-conformal mapping $\chi$ does not change the increase of the argument, when passing from $\nabla u$ to $\nabla h$. More precisely,   the chain rule yields
  \[
\partial_{z} u(z)= [\partial_w h](\chi(z)) [\partial_{z}  \chi](z),
\]
where $w$ is the space variable in $B(0,1)$. Thus, by using the properties of the logarithm,
\[
-\arg (\nabla u)=-\arg [(\nabla h)\circ \chi]+\arg(\partial_z \chi).
\]
But our assumption on $u$ guarantees that $[\partial_w h]\circ \chi$ is different from zero in a neighbourhood of $\partial D$, therefore $\partial_{z}  \chi$
is continuous on $\partial D$, does not vanish   on $\partial D$
and varies continuously with orientation preserved. This implies that
$\partial_z \chi$ returns to its initial value after a full turn, or equivalently
\[
\int_{+\partial D} d\arg(\partial_z \chi)=0.
\]
  These two last  identities  and  \eqref{defIu} prove \eqref{defIulocal}.}
  \end{proof}

\bt\label{t:nocriticalpoint}
Assume that the solution $u$ of
\eqref{EDPinu}-\eqref{bcu} 
belongs to $C^1(\bar \Omega)$
and satisfies \eqref{assumderiveenormaleextbdy} and 
 \eqref{assumderiveenormaleintbdy}. Then, recalling that  $M$ is the number of simply connected components of $\madeline{\mathcal{B}}$,
$u$ has $M-1$ critical points in $\Omega$ (when counted according to their multiplicities).
\et
\begin{proof}
By Lemmas \ref{l:ptcritiqueregularcase} and  \ref{l:ptcritiquenonregularcase},
we fix $\varepsilon>0$ small enough such that  {the set}
\[
\Omega_\varepsilon=\{x\in \Omega: \varepsilon<u(x)<1-\varepsilon\},
\]
is such that $\nabla u\ne \madeline{(0,0)^\top}$ in a neighborhood of $\partial\Omega_\varepsilon$.
Hence as $u\in C^1(\bar \Omega)$, by the implicit function theorem, $\Omega_\varepsilon$ has a $C^1$ boundary. Further again for $\varepsilon>0$ small enough, 
due to \eqref{assumderiveenormaleextbdy} and 
 \eqref{assumderiveenormaleintbdy}, we will have
\beq
 \label{assumderiveenormaleextbdyeps}
&&\partial_n u< 0 \hbox{ on } \Gamma_\varepsilon^{\rm ext} :=\{x\in \Omega: u(x)= \varepsilon\},
\\
&&\partial_n u>0 \hbox{ on } \Gamma_\varepsilon^{\rm int} :=\{x\in \Omega: u(x)= 1-\varepsilon\}.
\label{assumderiveenormaleintbdyeps}
\eeq
Now 
since the result of \cite{AlessandriniMagnanini:94} is based on oriented contour integrals (cfr Lemma \ref{l:indexu}), we proceed locally on simply connected oriented and disjoint patches $D_k$, $k=1,\cdots, K$
of $\Omega_\varepsilon$, whose union is $\Omega_\varepsilon$. Since the critical points of $u$ are finite, without loss of generality, we can further assume that
$\partial D_k\cap \partial D_{k'}$, with $k\ne k'$, has no critical point of $u$.
Further for each $k=1,\cdots, K$, we can assume that each $D_k$ is included into a larger simply connected open set
$\Omega_k'$ included into $\Omega$.
Therefore   the total number $N$ of critical points of $u$ in 
$\Omega$ is equal to the total number  of critical points of $u$ in 
$\Omega_\varepsilon$ and by  Lemma \ref{l:indexu},
$N$ is given by
\[
N=\sum_{k=1}^K I(D_k,u).
\]
Since on this sum, the contribution on $\partial D_k\cap \partial D_{k'}$ cancels, we get 
\[
N=-\frac{1}{2\pi}\int_{\madeline{+}\partial \Omega_\varepsilon} \madeline{d}\arg (\nabla u).
\]
Finally using  
 \eqref{assumderiveenormaleextbdyeps}, \eqref{assumderiveenormaleintbdyeps}
 and the boundary conditions on $u$ on $\Gamma_\varepsilon^{\rm ext}$ and $\Gamma_\varepsilon^{\rm int}$, we have that
 \beqs
 \int_{\madeline{+}\Gamma_\varepsilon^{\rm ext}} d\arg (\nabla u)=-1,
 \\
  \int_{\madeline{+}\Gamma_\varepsilon^{\rm int}} d\arg (\nabla u)=M.
\eeqs
Recall that $M$ is the number of connected components of $\madeline{\mathcal{B}}$.
This  implies that $N=M-1$.
\end{proof}


\section{Degenerate density in case of piecewise $C^{1,1}$ exterior boundary \label{s:degeneratedcase}}

In this section, we assume that {\bf (A1)} holds  with $I\geq 1$.
As explained in the introduction, in that case, $\rho(x)$ should tend to zero 
as $x$ approaches the corner $v_i$. Therefore we    are weakening the assumption
\eqref{bdrho} in the following way
\be\label{bdrhoweaker}
\forall \varepsilon>0, \exists \delta_\varepsilon>0:
\delta_\varepsilon \leq \rho \hbox{ a.e. in } \Omega\setminus \cup_{i=1}^I \bar B(v_i, \varepsilon).
\ee

In that situation, \eqref{bounds on u} does not follow from the maximum principle in $\Omega$, but it can be restored with the sole assumption that $u$ is positive in $\Omega$, as the next Lemma shows.

\bl\label{l:boundsonu}
Assume that the solution $u$ of
\eqref{EDPinu}-\eqref{bcu} 
belongs to $C^1(\bar \Omega)$ and that
\be\label{assumptionu>0}
u>0 \hbox{ in } \Omega.
\ee
Then \eqref{bounds on u}  remains valid.
\el
\begin{proof}
Since $u$ is continuous in $\bar \Omega$ and is zero on $\partial \mathcal{O}$,
given $\kappa\in (0,1)$,
there exists $\varepsilon_1>0$ small enough such that
\[
0\leq u(x)\leq \kappa, \forall x\in \cup_{i=1}^I (\bar B(v_i, \varepsilon_1)\cap \bar \Omega).
\]
Hence we can apply the maximum principle (see \cite[Theorems 8.1 and 8.19]{gilbarg:77}
for instance) to $u$ in $\Omega_{\varepsilon_1}:=\Omega\setminus \cup_{i=1}^I \bar B(v_i, \varepsilon_1)$ to find
\[
0<u<1 \hbox{ in } \Omega_{\varepsilon_1}.
\]
\madeline{Then} the two above properties yield $u<1$ in $\Omega$.
\end{proof}

\br{\rm
For the applications to system  \eqref{sys : gen}, the assumption  \eqref{assumptionu>0}
is natural. Indeed
if $\varphi\in C^1(\bar \Omega)$ is solution of \eqref{sys : gen} and \eqref{sys : DirBC},
then applying \cite[Theorem  2.2]{gilbarg:77}
to the first equation of  \eqref{sys : gen}, we find that $\varphi>0$ in $\Omega$.
}
\er

\subsection{Critical points are isolated}
Now the statement of Lemma \ref{l:ptcritiquenonregularcase} slightly differs 
and allows accumulation points near the corners.

\bl\label{l:ptcritiquenonregularcasedegeneraterho}
Let the assumptions {\bf (A1)}  with $I\geq 1$  and  \eqref{bdrhoweaker}  hold.
Assume that the solution $u$ of
\eqref{EDPinu}-\eqref{bcu} 
belongs to $C^1(\bar \Omega)$ 
and satisfies   \eqref{assumderiveenormaleextbdy},  
 \eqref{assumderiveenormaleintbdy}, and \eqref{assumptionu>0}. Then 
  $u$ has isolated critical points in $\bar \Omega$
  that can accumulate near the corners.
\el
\begin{proof}
It suffices to apply the arguments of the proof of Lemma \ref{l:ptcritiqueregularcase}
in $\Omega\setminus \cup_{i=1}^I \bar B(v_i, \varepsilon)$, for all $\varepsilon>0$.
\end{proof}

\subsection{Number of critical points inside the domain}

Since  we cannot directly use the argument principle\madeline{,} we first combine the use of the level lines and the maximum principle\madeline{. This allows us} to show that the number of critical points in $\Omega$
is finite and \madeline{hence} use the argument principle.

If $z_0\in \Omega$ is a critical point of index $I(z_0,u)$ of  a solution $u\in C^1(\bar \Omega)$  of
\eqref{EDPinu}-\eqref{bcu}, we recall from \cite[Lemma 2.5]{AlessandriniMagnanini:94}
that there exists a neighborhood $U\subset \Omega$ of $z_0$ such that the level line 
\[
L(z_0,u):=\{x\in \Omega: u(x)=u(\madeline{z_0})\},
\]
 is made of $I(z_0,u)+1$ simple arcs, whose pairwise intersection consists of $\{z_0\}$.
 Note further that in this neighborhood, $u-u(z_0)$ changes its sign when crossing one half-arc (due to \cite[\S 2.1]{Magnanini:16} and \cite[Theorem 2.1]{AlessandriniMagnanini:94}).
Then if $I(z_0,u)\geq 1$,  the sub-level set 
\[
 S_-(z_0,u)=\{x\in \Omega: u(x)< u(z_0)\},
 \]
and  the  super-level set
 \[
 S_+(z_0,u):=\{x\in \Omega: u(x)> u(z_0)\}
 \]
can be  decomposed into different open, connected, disjoint and non empty subsets,
their number being possibly infinite.
Let us show that this is not the case and further give some specific properties of these connected sets.

\bl\label{l:propS_k}
Let $z_0\in \Omega$ be a critical point of index $I(z_0,u)\geq 1$ of  {a solution} $u\in C^1(\bar \Omega)$ of
\eqref{EDPinu}-\eqref{bcu}.
Then  {we have that}
\begin{enumerate}
\item
$S_-(z_0,u)$ is connected, and its boundary contains the boundary of $\mathcal{O}.$
\\
\item
the set $S_+(z_0,u)$ is made of a finite number $K_+(z_0,u)$ of connected and disjoint non empty subsets,
 \[
S_+(z_0,u)=\cup_{k=1}^{K_+(z_0,u)}S_{k,+}(z_0,u),
\]
and for all $k=1, \cdots, K_+(z_0,u)$, the boundary of
$
S_{k,+}(z_0,u)$ contains the boundary of at least one $B_m$. 
\end{enumerate}
Consequently, 
\be\label{bdK}
K_+(z_0,u)\leq M.
\ee
Further we have that
\be\label{bdKI}
K_+(z_0,u)\geq I(z_0,u)+1,
\ee
%
 and therefore 
 \be\label{bdI}
 I(z_0,u)\leq M-1.
 \ee
\el
\begin{proof}
As $z_0\in \Omega$,  {then}
$u(z_0)\in (0,1)$ due to  \eqref{bounds on u}, and by the boundary condition  \eqref{bcu}, we get
\be\label{proplevelline}
L(z_0,u) \subset \Omega.
\ee
This directly implies that 
for any open, connected, and non empty  component $V_-$ of $S_{-}(z_0,u)$ (resp. $V_+$ of $S_{+}(z_0,u)$),  $\overline{V_\pm}$ cannot be included into $\Omega$. Indeed, if 
$\overline{V_\pm}\subset \Omega$, then its boundary  would be  included into $L(z_0,u)$.
Hence using the maximum principle on $u$ in $\overline{V_\pm}$, we deduce that
$u=u(\madeline{z_0})$ in $\madeline{V_\pm}$, which contradicts Lemma \ref{l:ptcritiquenonregularcasedegeneraterho}.
The assertions concerning the sets $S_-(z_0,u)$ and  $S_{+}(z_0,u)$ and the bound on $K_+(z_0,u)$ directly follow from these properties. 
Now \eqref{bdKI} directly follows from the behavior of  $u-u(z_0)$  in a neighborhood of $z_0$ described previously, because in that neighborhood  there are $I(z_0,u)+1$ open disjoint sets where 
$u-u(z_0)>0$.
\end{proof}

%

\bt\label{t:nocriticalpointdegenerated}
Let the assumptions {\bf (A1)} with $I\geq 1$  and  \eqref{bdrhoweaker}  hold.
Assume that the solution $u$ of
\eqref{EDPinu}-\eqref{bcu} 
belongs to $C^1(\bar \Omega)$ 
and satisfies  \eqref{assumderiveenormaleextbdy},
 \eqref{assumderiveenormaleintbdy}, and \eqref{assumptionu>0}. Then recalling that  $M$ is the number of simply connected components of $\madeline{\mathcal{B}}$,
$u$ has $M-1$ critical points in $\Omega$ (when counted according to their multiplicities).
\et
\begin{proof}
1. If $\madeline{\mathcal{B}}$ is simply connected,  then $M=1$,  {and hence}
 the bound \eqref{bdI} directly forbids the existence of $z_0\in \Omega$ such that  $I(z_0,u)\geq 1$.
\\
2. Assume now that $M\geq 2$, and that  $u$ has a critical point $z_0\in \Omega$ of index $I(z_0,u)\geq 1$. 
Then we
consider the set of all critical points $z_j\in \Omega$ of index $I(z_j,u)\geq 1$ 
such that $u(z_j)=u(z_0)$, $j=1,\cdots, J$, with $J\in \N$. We use  the convention that $J=0$
when such  critical points do not exist. 

Now assume that there exists another  critical point $z^*\in \Omega$ of index $I(z^*,u)\geq  1$. Then we can distinguish between the following two cases:
\begin{enumerate}[(i)]
\item
 if $u(z^*)<u(z_0)$, then $L(z^*,u)$ is contained in $S_-(z_0, u)$. Therefore, 
$S_-(z^*, u) \subset S_-(z_0, u)$ and
the boundary of one $S_{k,+}(z^*, u)$, with $k\geq 1$, say $S_{1,+}(z^*, u)$, contains $\partial \madeline{\mathcal{B}}$, otherwise $L(z^*,u)$ would cross 
$L(z_0,u)$. This means that the boundary of $S_{2,+}(z^*, u)$ is contained in $L(z^*,u)$,
and this is impossible, due to the maximum principle.

\item
if  $u(z^*)>u(z_0)$ then, on one hand, $L(z^*,u)$ is contained in one $S_{k,+}(\madeline{z_0}, u)$, with $k\geq 1$ and, on the other hand, since the boundary of $S_-(z^*,u)$ contains $\partial\mathcal{O}$, 
$L(z^*,u)$ has to cross $L(z_0,u)$. This is a contradiction.
\end{enumerate}

In conclusion such a critical point $z^*$ does not exist, and therefore
$u$ has a finite number of critical points in $\Omega$. This implies that we can now apply the arguments of the proof of Theorem \ref{t:nocriticalpoint}, by simply taking the 
simply connected open sets
$\Omega_k'$ contained in
$\Omega\setminus \cup_{i=1}^I \bar B(v_i, \varepsilon')$, with $\varepsilon'>0$ small enough.
\end{proof}

\section{Checking the Hopf lemma on the exterior/interior boundary\label{s:Hopf}}

{As already mentioned in the introduction, we recall that the Hopf lemma may not hold when $\rho$ is not H\"older continuous. 
Thus, in this section we  give some sufficient conditions on the regularity of the boundary of $\mathcal{O}$ (resp. $\mathcal{B}$)
that guarantee that \eqref{assumderiveenormaleextbdy} (resp. \eqref{assumderiveenormaleintbdy}) holds. More precisely, we suppose that the boundary of $\mathcal{O}$ is piecewise analytic, except at a finite number
 of points, where the interior angle formed by the two tangents is equal to $\frac{\pi}{2}$. Namely, we assume the following two conditions and (A1).}

\medskip\noindent
{\bf (A2)}  There exist
$\madeline{a_i}\in \partial \mathcal{O}\setminus \cup_{\madeline{j}=1}^I \{v_j\}$, $i=I+1,\cdots, I+J$, with $J\in \N^*$, where the boundary  coincides locally with the image by a conformal map of  a half plane (identifying $\R^2$ with $\C$). Namely,  there exist
a neighbourhood $U_i$ of \madeline{$a_i$} in $\C$,
an open, non empty,  simply connected and bounded set $V_i$ of $\C$ such that $0\in V_i$
and a conformal 
mapping 
\[
\Phi_i:V_i\cap \madeline{\sigma_+}\to U_i\cap \Omega,
\]
 where \madeline{$\sigma_+ :=\{(x_1,y_1)\in \R^2: x_1>0\} $} is the right-hand half-plane.

\medskip\noindent
{\bf (A3)} 
 The set $\cup_{i=1}^{I+J} \overline{U_i\cap \Omega}$ is covering the whole
 boundary of $\mathcal{O}$.
 
 \bt[Hopf lemma on the exterior boundary]\label{t:Hopfextbdy}
{ Assume  
that conditions {\bf (A2)} to  {\bf (A3)} are satisfied, together with {\bf (A1)}. Suppose that $\rho\in L^\infty(\Om)$ and assume that \eqref{bdrho} holds.}
 Let $u\in C^1(\bar \Omega)$ be   the solution of
\eqref{EDPinu}-\eqref{bcu}.
Then 
\eqref{assumderiveenormaleextbdy} holds.
\et
\begin{proof}
1. Let us start by proving \eqref{assumderiveenormaleextbdy} when {\bf (A2)} to  {\bf (A3)} are satisfied and $\partial \mathcal{O}$ is $C^{1,1}$. 
\\
a. For a point $\madeline{a_i}$, with $i\ge 1$, we  first straighten the boundary of $U_i\cap \mathcal{O}$  by means of  the conformal map $\madeline{\Phi_i}$. This  {procedure} directly leads to problem  \eqref{pb**}.
 Hence an odd reflection  with respect to the axis $\re(w)=0$ leads to the solution
 $u^{***}\in H^1(V^{***})$  of
\eqref{EDP***} that satisfies  
 $u^{***}=0$ on the segment $S_1:=\overline{V^{***}}\cap \{ x_1=0\}$ 
 and 
 \[
0<|u^{***}|<1 \hbox{ on } \overline{ V^{***}}\setminus S_1.
\]
Hence by \cite[Lemma 2.5]{AlessandriniMagnanini:94} we conclude that 
$\nabla u^{***}\ne (0,0)^\top$ along $S_1$ and consequently $\nabla u\ne (0,0)^\top$ along $U_i\cap \partial\mathcal{O}$.\\
b. Conclusion: By the assumption {\bf (A3)}, 
\[
\nabla u\ne (0,0)^\top \hbox{ on } \partial \mathcal{O},
\]
and as the second boundary condition from \eqref{bcu} implies that the tangential derivative is zero on $\partial \mathcal{O}$, we can conclude that 
\eqref{assumderiveenormaleextbdy} holds.
\\
2. We now prove \eqref{assumderiveenormaleextbdy} when {\bf (A1)} with $I\geq 1$, {\bf (A2)} and {\bf (A3)} are satisfied.
\\
a. For a point $\madeline{a_i}$, with $i\ge I+1$, the arguments of point 1 yield $\nabla u\ne (0,0)^\top$ along $U_i\cap \partial\mathcal{O}$.
\\
b. For a corner $v_i$, $i=1,\cdots, I$, we follow the arguments of the proof of Lemma \ref{l:ptcritiquenonregularcase}. More precisely, with the different changes of variables, we arrive  {at showing that $u^{***}\in H^1(V^{***})$ and is  a} solution of
\eqref{EDP***} that satisfies  {\eqref{eq:12} and that}
 $u^{***}=0$ on the segments $S_1:=\overline{V^{***}}\cap \{ x_1=0\}$ and $S_2=\overline{V^{***}}\cap \{ y_1=0\}$. This implies that $0$ is a critical point of $u^{***}$
 and by \cite[Lemma 2.5]{AlessandriniMagnanini:94} applied to $u^{***}$ in $V^{***}$, the points of $S_1\setminus\{0\}$ and $S_2\setminus\{0\}$ cannot be critical points of index $\geq 1$. Hence $\nabla u^{***}\ne (0,0)^\top$ along $S_1\setminus\{0\}$ and $S_2\setminus\{0\}$.
 Coming back to $u$, we get that $\nabla u\ne (0,0)^\top$ along $(U_i\cap \partial\mathcal{O})\setminus\{v_i\}$.
 \\
c. Conclusion: By the assumption {\bf (A3)}, 
\[
\nabla u\ne (0,0)^\top \hbox{ on } \partial \mathcal{O}\setminus \cup_{i=1}^I \{v_i\},
\]
and as the second boundary condition from \eqref{bcu} implies that the tangential derivative is zero on $\partial \mathcal{O}\setminus \cup_{i=1}^I \{v_i\}$, we can conclude that 
\eqref{assumderiveenormaleextbdy} holds.
\end{proof}

This result actually extends to the case of a degenerate density.
 
 \bc\label{c:Hopfextbdy}
 Let the assumptions  {\bf (A1)} with $I\geq 1$, {\bf (A2)}  and  {\bf (A3)} be satisfied and assume that \eqref{bdrhoweaker} holds.
 Let   $u\in C^1(\bar \Omega)$ be   the solution of
\eqref{EDPinu}-\eqref{bcu} that satisfies \eqref{assumptionu>0}.
Then 
\eqref{assumderiveenormaleextbdy} holds.
\ec
\begin{proof}
Only step 1 of the previous Theorem has to be adapted. Indeed in the first step, the sole difference is that
the density $\rho^{***}$ in
\eqref{EDP***} degenerates near $0$. But due to  \eqref{bdrhoweaker}, we can apply 
\cite[Lemma 2.5]{AlessandriniMagnanini:94}  in the simply connected domains
\[
 V^{***}\cap \{(x_1,y_1)\in \R^2: y_1>-|x_1]\}\setminus B(0,\varepsilon),
 \]
 and 
 \[
 V^{***}\cap \{(x_1,y_1)\in \R^2: y_1<|x_1]\}\setminus B(0,\varepsilon),
 \]
with $\varepsilon>0$ as small as we want, to deduce that 
$\nabla u^{***}\ne (0,0)^\top$ along $S_1\setminus B(0,\varepsilon)$ and $S_2\setminus B(0,\varepsilon)$. As $\varepsilon\to 0^+$, we deduce as before that 
$\nabla u^{***}\ne (0,0)^\top$ along $S_1\setminus\{0\}$ and $S_2\setminus\{0\}$.
\end{proof}

Obviously adapting the assumptions {\bf (A2)} to  {\bf (A3)} 
to $\partial \mathcal{B}$, we directly get \eqref{assumderiveenormaleintbdy}.
More precisely, we require:

\medskip\noindent
{\bf (A4)} there exist
$v^*_k\in \partial \madeline{\mathcal{B}}$, $k=1,\cdots, J$, with $K\in \N^*$, where the boundary  coincides locally with the image by a conformal map of  a half plane (identifying $\R^2$ with $\C$), namely there exist
a neighbourhood $U^*_k$ of $v^*_k$ in $\C$,
an open, non empty,  simply connected and bounded set $V^*_k$ of $\C$ such that $0\in V^*_k$
and a conformal 
mapping 
\[
\Phi^*_k:V^*_k\cap \pi_+\to U_k^*\cap \Omega.
\]
 
\medskip\noindent
 {\bf (A5)}  The set $\cup_{k=1}^{K} \overline{U^*_k\cap  \Omega}$ is covering the whole
 boundary of $\madeline{\mathcal{B}}$.

\bc[Hopf lemma on the interior boundary]\label{c:Hopfintbdy}
 Let the assumptions  {\bf (A4)} to  {\bf (A5)} be satisfied and assume that \eqref{bdrhoweaker}holds.
 Let  $u\in C^1(\bar \Omega)$ be   the solution of
\eqref{EDPinu}-\eqref{bcu} that satisfies \eqref{assumptionu>0}.
Then 
\eqref{assumderiveenormaleintbdy} holds.
\ec
\begin{proof}
We apply the above arguments to $1-u$.
\end{proof}

\section{Numerical illustrations \label{sexample}}

In this section, we present different numerical tests that confirm our theoretical results.
To compute the numerical approximation of the solution of \eqref{EDPinu}-\eqref{bcu}
with different values of $\rho$, we use the Finite Element Method with  the help of the library \textit{Netgen} \cite{Netgen} from \textit{python}. Here we restrict ourselves to  piecewise $\mathbb{P}_1$ finite elements   with a mesh of meshsize $h=0.005$.

\subsection{An annulus} \label{ss:annulus}
\begin{figure}[h]
\begin{minipage}[c]{.46\linewidth}
        \centering
        \includegraphics[scale=1.5,width=8cm]{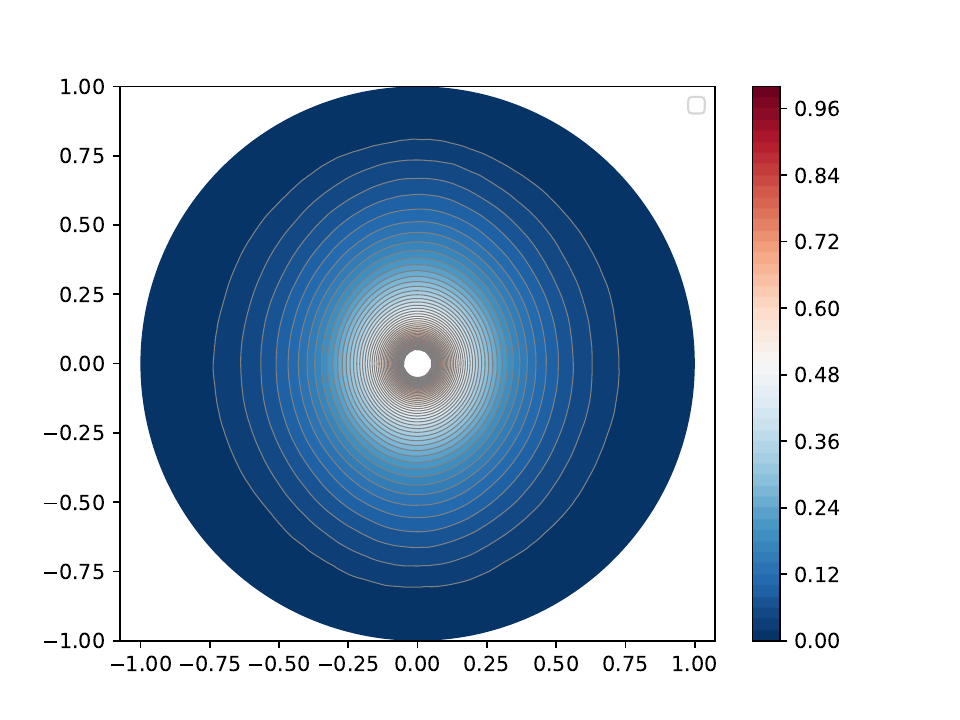}
    \end{minipage}
    \hfill%
    \begin{minipage}[c]{.46\linewidth}
        \centering
        \includegraphics[scale=1.5,width=8cm]{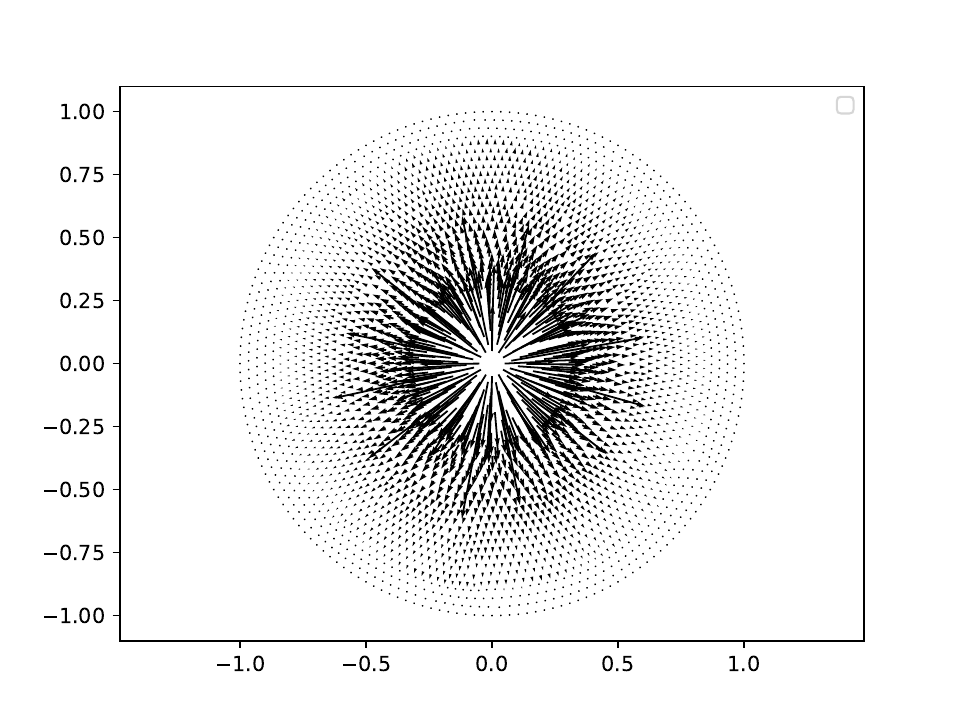}
    \end{minipage}
     \caption{An annulus and a smooth $\rho$: 20 level lines (left) and a gradient map (right)\label{anneau_x2plus1}}
\end{figure}

We first start with the case of an annulus $\Omega$ of center 0, of interior radius $0.05$
and exterior radius 1, where $\Gamma_c$ corresponds to the boundary of the disc
$D(0, 0.05)$ and $\Gamma_d$ corresponds to the boundary of the disc
$D(0, 1)$. We first take for $\rho$ the smooth function and uniformly bounded from below:
\[
\rho(x_1,x_2)= x_1^2 +1/8.
\]
In this situation, the solution $u$ of
\eqref{EDPinu}-\eqref{bcu}  clearly 
belongs to $C^1(\bar \Omega)$ by standard regularity results. Further as mentionned in the introduction as $\rho\in C^{0,\alpha}(\bar \Omega)$, the assumptions
\eqref{assumderiveenormaleextbdy},
 \eqref{assumderiveenormaleintbdy} are satisfied.
Therefore  by Theorem 
\ref{t:nocriticalpoint}, $u$ has no critical points in $\Om$.
We then   compute an approximation of  this solution using the numerical scheme described above. 
In Figure \ref{anneau_x2plus1},
we have shown 
20 level lines of $u$ and a gradient map on a coarse mesh. Both figures confirm that $u$ has indeed no critical points in $\Omega$.

%
%

\begin{figure}[h]
\begin{minipage}[c]{.46\linewidth}
        \centering
        \includegraphics[scale=1.5,width=8cm]{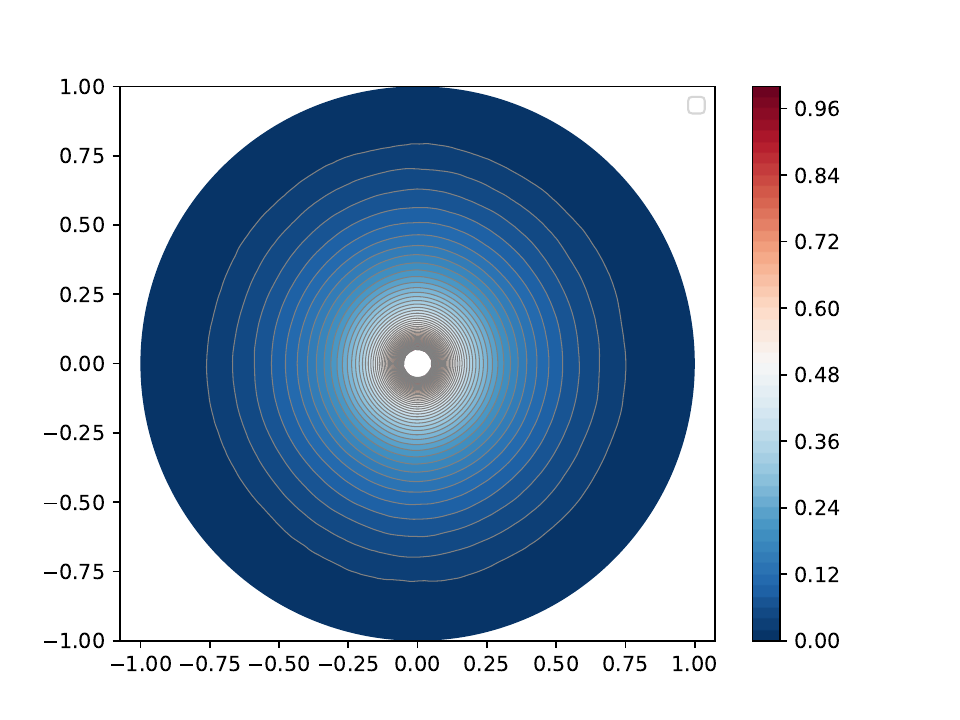}
    \end{minipage}
    \hfill%
    \begin{minipage}[c]{.46\linewidth}
        \centering
      \includegraphics[scale=1.5,width=8cm]{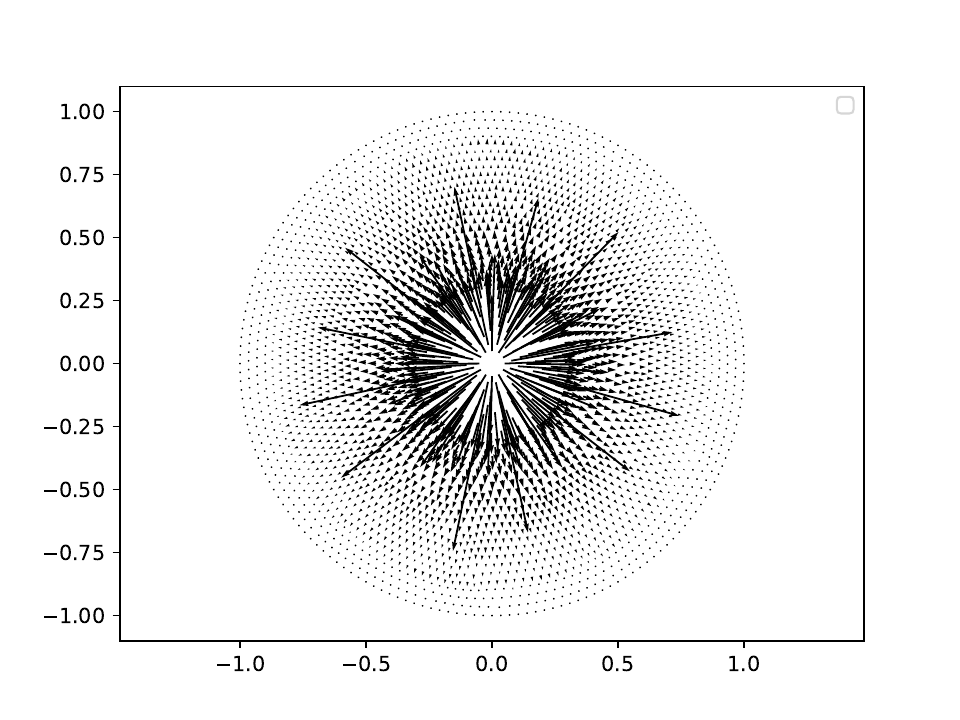}
    \end{minipage}
     \caption{An annulus and a nonsmooth  $\rho$: 20 level lines (left) and a gradient map (right)\label{anneau_absxplus1}}
\end{figure}

In  Figure \ref{anneau_absxplus1}, we can see that the same phenomenon occurs for the same domain $\Omega$ and  the choice of 
\[
\rho(x_1,x_2)= |x_1|+1/8,
\]
that is only Lipschitz, but is still uniformly bounded from below. With such a choice using odd and even reflections with respect to the $x_1$ axis, we can show that the solution $u$ of   
\eqref{EDPinu}-\eqref{bcu}  
belongs to $C^1(\bar \Omega)$, while it satisfies \eqref{assumderiveenormaleextbdy},
 \eqref{assumderiveenormaleintbdy} due to the Lipschitz property of $\rho$.
%
%
%

 \subsection{The unit disc with three holes}
\begin{figure}[h]
\begin{minipage}[c]{.46\linewidth}
        \centering
        \includegraphics[scale=1.5,width=8cm]{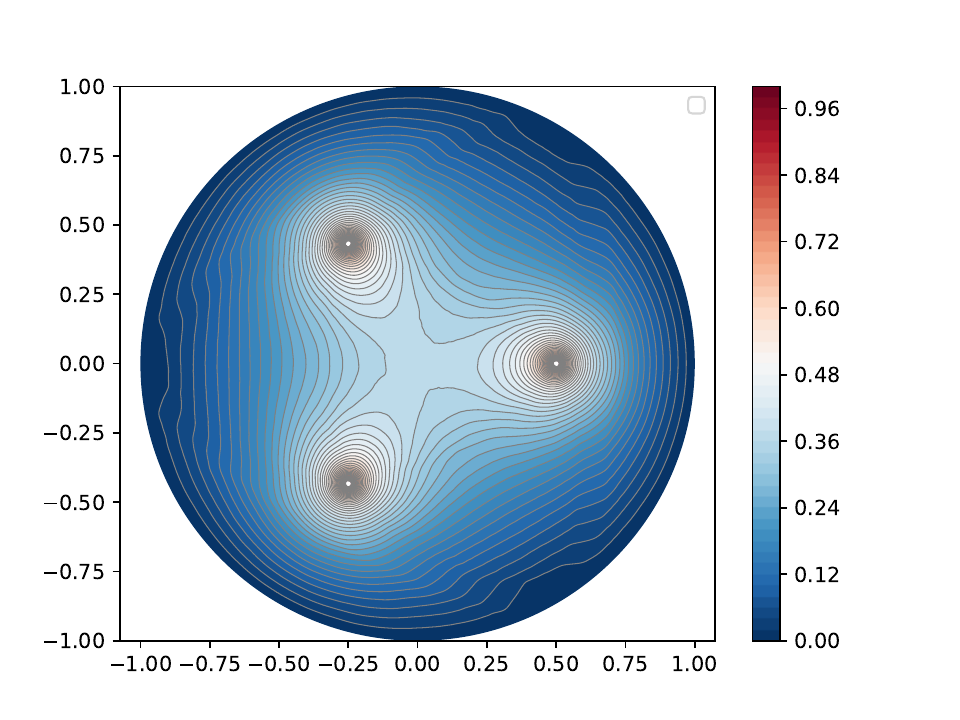}
    \end{minipage}
    \hfill%
    \begin{minipage}[c]{.46\linewidth}
        \centering
       \includegraphics[scale=1.5,width=8cm]{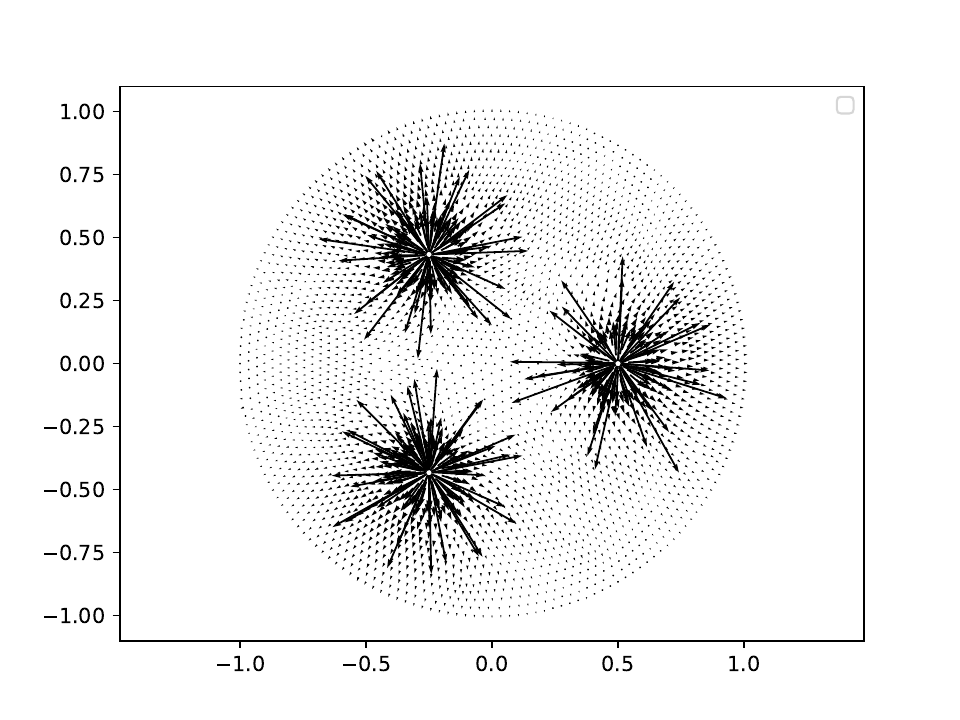}
    \end{minipage}
     \caption{The unit ball with 3 holes and a smooth  $\rho$: 20 level lines (left) and a gradient map (right)\label{anneau_x2plusy2plus1}}
\end{figure}

 Now we consider the domain $\Omega=B(0,1)\setminus \cup_{i=1}^3 \bar B(c_i, 0.01)$, where 
 \[
 c_1=(0.5,0), c_2=0.5(\cos(2\pi/3),\sin(2\pi/3)),  c_3=0.5(\cos(4\pi/3),\sin(4\pi/3)),
 \] and take first 
  the smooth function and uniformly bounded from below:
\be\label{rhpr2+1}
\rho(x_1,x_2)= x_1^2+ x_2^2 +1.
\ee
Clearly the assumptions on the solution $u$ of  
\eqref{EDPinu}-\eqref{bcu}  in  Theorem 
\ref{t:nocriticalpoint} are satisfied and therefore $u$ has two critical points in $\bar \Om$. But the choice of $\rho$ and of the holes imply that this solution is invariant by a rotation of $2\pi/3$, therefore $(0,0)$ is the unique critical point of geometric index 2.
Again, we have   computed an approximation of  this solution and 20 level lines and a gradient map on a coarse mesh of $u$ are presented in Figure \ref{anneau_x2plusy2plus1}
and suggest that    indeed the  critical point of  $u$ in $\Omega$  is $(0,0)$.

%
%

\begin{figure}[h]
\begin{minipage}[c]{.46\linewidth}
        \centering
       \includegraphics[scale=1.5,width=8cm]{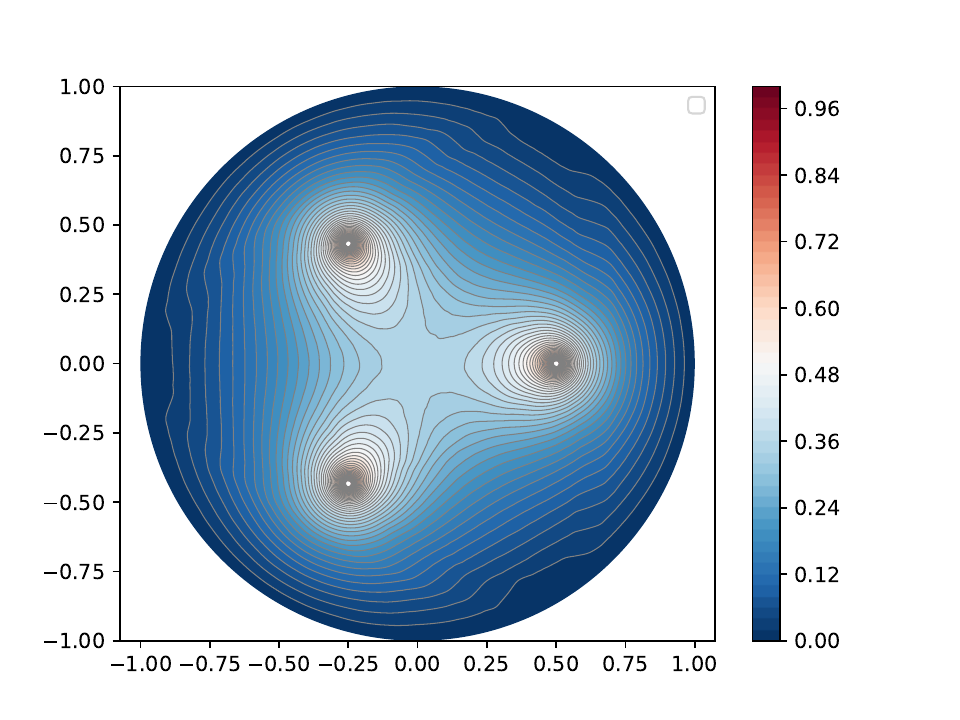}
    \end{minipage}
    \hfill%
    \begin{minipage}[c]{.46\linewidth}
        \centering
       \includegraphics[scale=1.5,width=8cm]{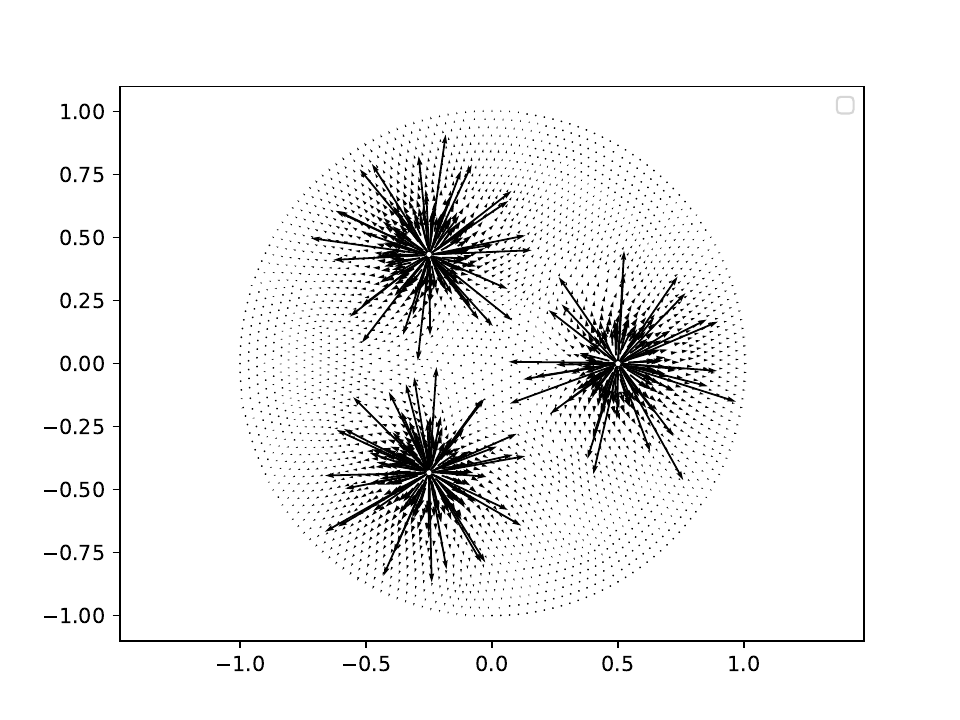}
    \end{minipage}
     \caption{The unit ball with 3 holes and a degenerate   $\rho$: 20 level lines (left) and a gradient map (right)\label{anneau_x2plusy2}}
\end{figure}

For the same configuration, we have also tested the case 
\[
\rho(x_1,x_2)= \sqrt{x_1^2+ x_2^2},
\]
that is not covered by our theory because $\rho$ degenerates at an interior point of the domain. Nevertheless a solution of problem
\eqref{EDPinu}-\eqref{bcu}
exists in some weighted Sobolev space, see \cite{FabesKenigSerapioni:82,DeCiccoVivaldi:96} for instance, and
 Figure \ref{anneau_x2plusy2}  suggests that     the unique critical point of  $u$ in $\Omega$  is $(0,0)$.

%
%

 \subsection{A half-disc with one hole}
\begin{figure}[h]
\begin{minipage}[c]{.46\linewidth}
        \centering
         \includegraphics[scale=1.5,width=8cm]{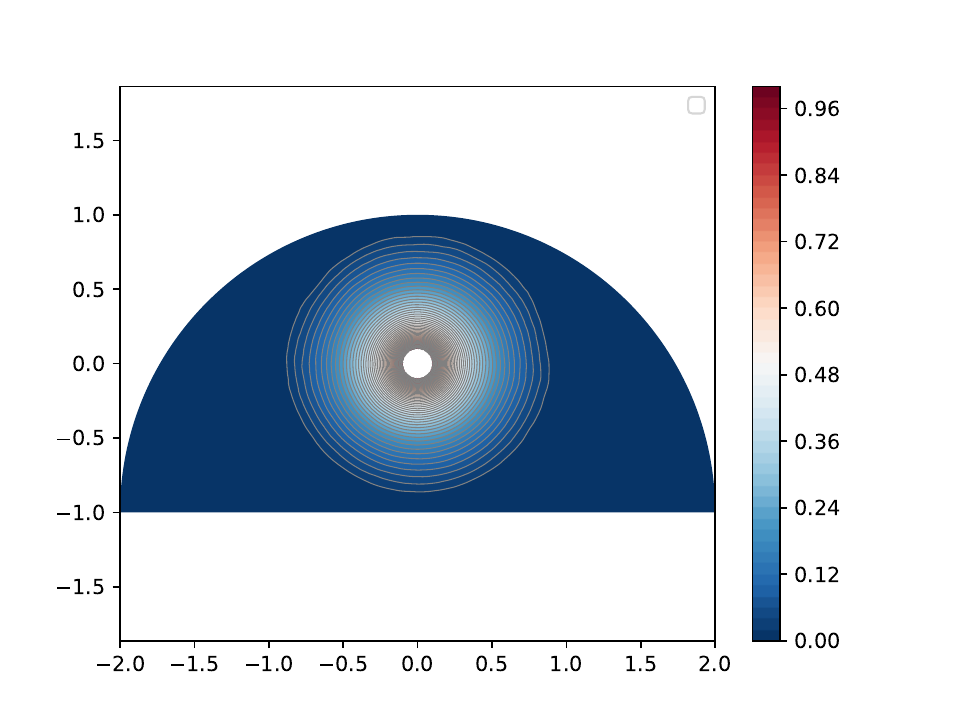}
    \end{minipage}
    \hfill%
    \begin{minipage}[c]{.46\linewidth}
        \centering
      \includegraphics[scale=1.5,width=8cm]{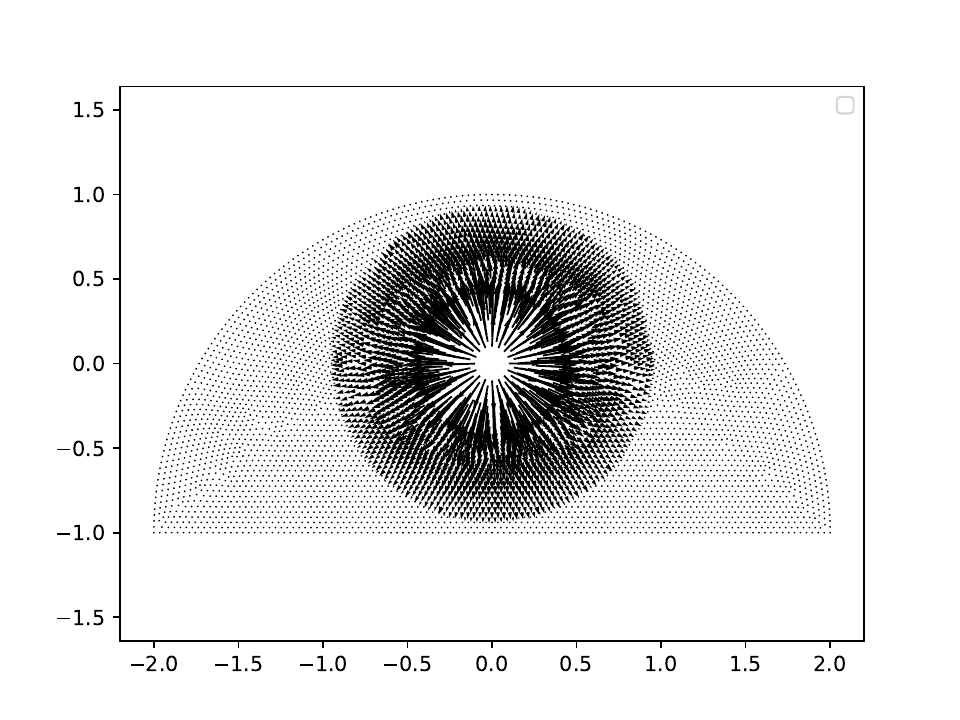}
    \end{minipage}
     \caption{A half-disc  with 1 hole and  a smooth  $\rho$: 20 level lines (left) and a gradient map (right)\label{halfdisconehole_x2plusy2+1}}
\end{figure}

 Now we consider the domain $\Omega=B_+\setminus  \bar B((0,0), 0.01)$, where 
 \be\label{halfball}
 B_+=\{(x_1,x_2)\in B((0,-1), 2): x_1>-1\},
 \ee
 and we  first take  $\rho$ defined by 
 \eqref{rhpr2+1}. Again as $\rho$ is smooth,  the assumptions on the solution $u$ of  
\eqref{EDPinu}-\eqref{bcu}  in  Theorem 
\ref{t:nocriticalpoint} are satisfied, and therefore $u$   has no critical points in $\Om$.
This property is confirmed by Figures \ref{halfdisconehole_x2plusy2+1}.

%
%

\begin{figure}[h]
\begin{minipage}[c]{.46\linewidth}
        \centering
      \includegraphics[scale=1.5,width=8cm]{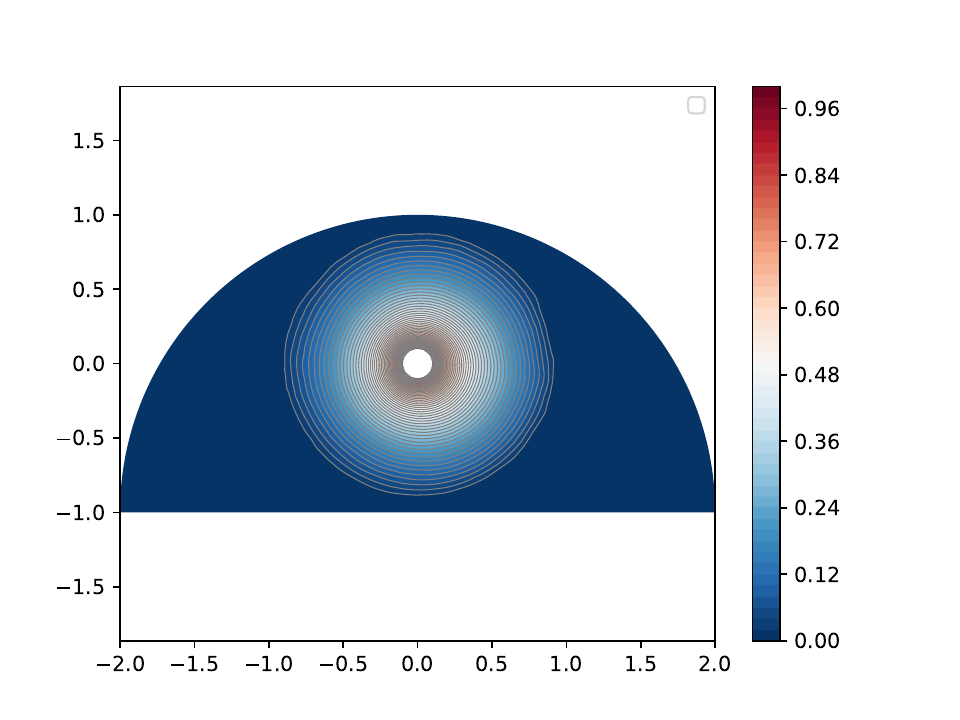}
    \end{minipage}
    \hfill%
    \begin{minipage}[c]{.46\linewidth}
        \centering
      \includegraphics[scale=1.5,width=8cm]{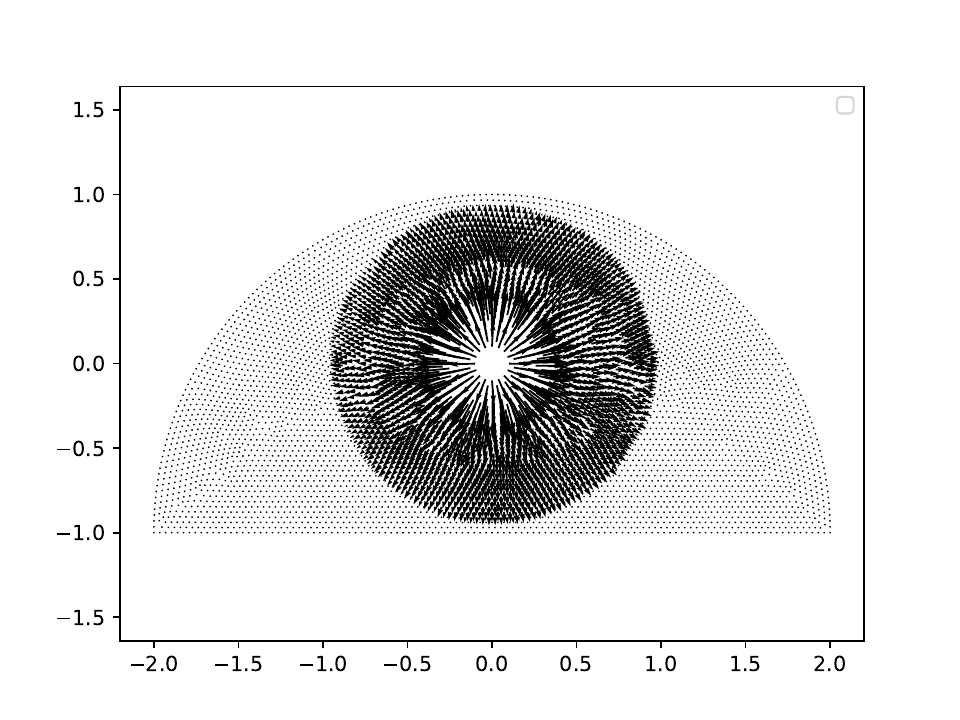}
    \end{minipage}
     \caption{A half-disc  with 1 hole and  a degenerate $\rho$: 20 level lines (left) and a gradient map (right)\label{halfdisconehole_xmoins2carre_plusyplus1carre}}
\end{figure}

Let us now take 
\[
\rho(x_1,x_2)=\sqrt{(x_1-2)^2+(x_2+1)^2},
\]
that   degenerates at the vertex $(2,-1)$ of $\Om$ but satisfies \eqref{bdrhoweaker}.
In such a case, a  solution  of problem
\eqref{EDPinu}-\eqref{bcu}
exists in some weighted Sobolev space, see again 
\cite{FabesKenigSerapioni:82,DeCiccoVivaldi:96} for instance, hence its $C^1$ regularity is not guaranteed. Nevertheless  Figure \ref{halfdisconehole_xmoins2carre_plusyplus1carre}
  seems to indicate that the solution  has no critical points in $\Om$.

%
%

 \subsection{An half-disc with three holes}

\begin{figure}[h]
\begin{minipage}[c]{.46\linewidth}
        \centering
      \includegraphics[scale=1.5,width=8cm]{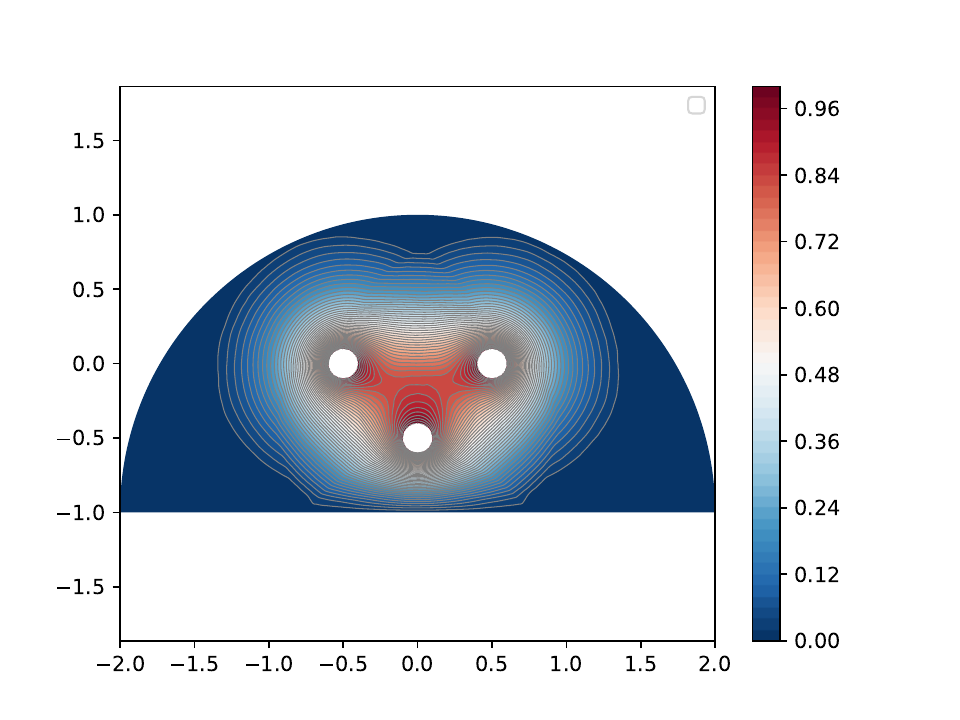}
    \end{minipage}
    \hfill%
    \begin{minipage}[c]{.46\linewidth}
        \centering
      \includegraphics[scale=1.5,width=8cm]{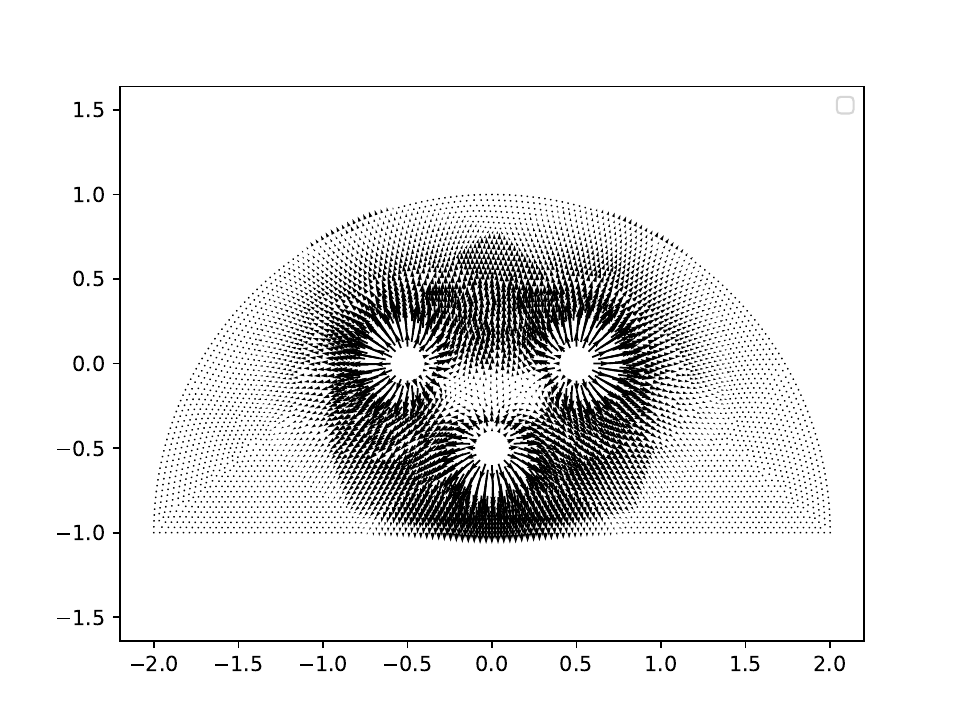}
    \end{minipage}
     \caption{A half-disc  with 3 holes and  a smooth $\rho$: 20 level lines (left) and a gradient map (right)\label{halfdisc3hole_x2plusy2+1}}
\end{figure}

 Now we consider the domain $\Omega=B_+\setminus \cup_{i=1}^3 \bar B(c_i, 0.01)$, where 
$B_+$  is given by  \eqref{halfball},   $c_1=(-0.5,0), c_2=(0, -0,5)$, and $c_3=(0.5,0)$ and take  $\rho$ defined by 
 \eqref{rhpr2+1}. Again as $\rho$ is smooth,  the assumptions on the solution $u$ of  
\eqref{EDPinu}-\eqref{bcu}  in  Theorem 
\ref{t:nocriticalpoint} are satisfied, and therefore $u$   has two critical points in $\Om$.
This property is confirmed by Figure \ref{halfdisc3hole_x2plusy2+1} where the   critical points seem to be near $(-0.25, -0,25)$ and $(0.25, -0,25)$.

%
%

 \subsection{An annulus with discontinuous coefficients} 
 
 In this last subsection, we want to treat the case where the coefficient $\rho$ is discontinuous. In this case, as mentionned in section \ref{s:setting},
 the solution of  
\eqref{EDPinu}-\eqref{bcu}  is in general no more in $C^1(\bar \Omega)$, see Remark \ref{rk:reg}. Therefore  even if we  are  slightly outside our regularity assumptions, we want to check  analytically or numerically if the solution of  
\eqref{EDPinu}-\eqref{bcu}  may have or not critical points for with a piececewise constant coefficient $\rho$.

For that purpose, we consider the annulus of subsection 
 \ref{ss:annulus} and take different choices of $\rho$.
 
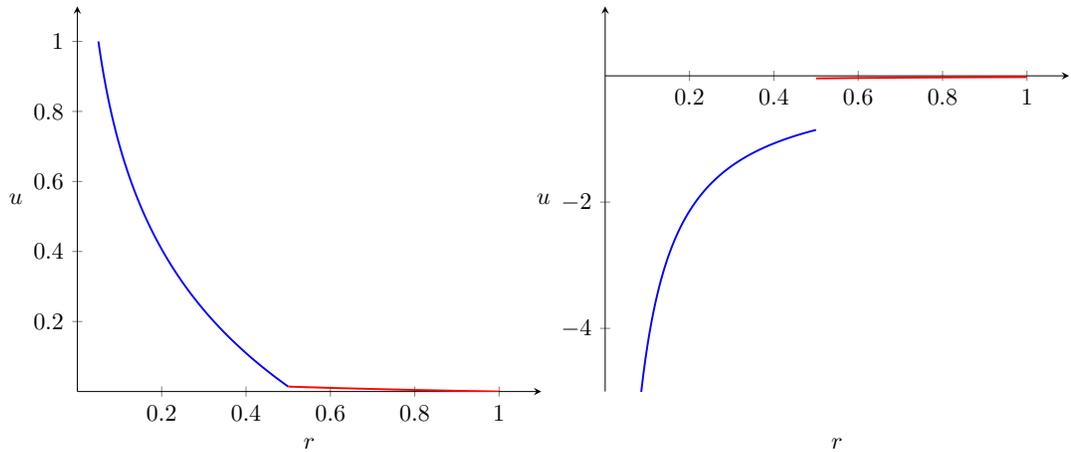
\begin{figure}[h!t]
\begin{minipage}[c]{.46\linewidth}
  \centering
  \begin{tikzpicture}[scale = 0.9]
    \begin{axis}[
      axis lines=middle,
      xlabel=$r$, ylabel=$u$,
      xmin=0, xmax=1.1,
      ymin=0, ymax=1.1,
      samples=100,
      domain=0.01:1.1,  
      xlabel style={anchor=north, at={(axis description cs:0.5,-0.1)}},
      ylabel style={anchor=east, at={(axis description cs:-0.1,0.5)}}    
    ]
    
    \addplot[blue, thick, domain=0.05:0.5] {1+ 1/(1/21 * ln(0.5)-ln(0.5/0.05)) * ln(x/0.05)};

     \addplot[red, thick, domain=0.5:1] {1/21 *1/(1/21 * ln(0.5)-ln(0.5/0.05)) * ln(x)};

  \end{axis}
  \end{tikzpicture}
  \end{minipage}
\begin{minipage}[c]{.46\linewidth}
  \centering
  \begin{tikzpicture}[scale = 0.9]
    \begin{axis}[
      axis lines=middle,
      xlabel=$r$, ylabel=$u$,
      xmin=0, xmax=1.1,
      ymin=-5, ymax=1.1,
      samples=100,
      domain=0.01:1.1,  
      xlabel style={anchor=north, at={(axis description cs:0.5,-0.1)}},
      ylabel style={anchor=east, at={(axis description cs:-0.1,0.5)}}    
    ]
    
    \addplot[blue, thick, domain=0.05:0.5] {1/(1/21 * ln(0.5)-ln(0.5/0.05)) * 1/ x};

    \addplot[red, thick, domain=0.5:1] {1/21 *1/(1/21 * ln(0.5)-ln(0.5/0.05)) * (1/ x)};
  \end{axis}
  \end{tikzpicture}
  \end{minipage}
  \caption{The solution $u$ (left) and  $\frac{\partial u}{\partial r}$ (right) with respect to $r$
  for a radially discontinuous coefficients\label{fig:madeline}}
\end{figure}

 The first choice consists in taking a radially symmetric coefficient $\rho$, namely
 \[
 \rho(x)=\left\{
 \begin{array}{ll}
 \rho_- & \text{ if } r_0<|x|<r_1,
 \\
 \rho_+& \text{ if } r_1<|x|<1,
 \end{array}
 \right.
\]
where $ \rho_-,  \rho_+$ are two positive real numbers and $0<r_0<r_1<1$ (here $r_0$ can be arbitrary). In that case, we readily check that the exact solution $u$ of \eqref{EDPinu}-\eqref{bcu} is radial and is given  in polar coordinates $(r,\theta)$ by
\[
u(x)=\left\{
 \begin{array}{ll}
1+a_- \ln \left(\frac{r}{r_0}\right)& \text{ if } r_0<r<r_1,
 \\
a_+ \ln r& \text{ if } r_1<r<1,
 \end{array}
 \right.
\]
where 
\[ a_+=\frac{\rho_-}{\rho_+} a_-\text{ and }
a_-=\left(\frac{\rho_-}{\rho_+} \ln r_1-\ln \left(\frac{r_1}{r_0}\right)\right)^{-1}
\]
are both negative. This implies that 
\[
\nabla u(x)=\left\{
 \begin{array}{ll}
\frac{a_-}{r} \left(
 \begin{array}{ll}
 \cos \theta\\
 \sin \theta
 \end{array}
\right)& \text{ if } r_0<r<r_1,
 \\
\frac{a_+}{r} \left(
 \begin{array}{ll}
 \cos \theta\\
 \sin \theta
 \end{array}
\right)& \text{ if } r_1<r<1,
 \end{array}
 \right.
\]
which means that $u$ has no critical points, see Figure \ref{fig:madeline} for an illustration with $r_0=0.05$, $r_1=0.5$, $\rho_-=1$, and   $\rho_+=21$.

\begin{figure}[h]
\begin{minipage}[c]{.46\linewidth}
        \centering
       \includegraphics[scale=1.5,width=8cm]{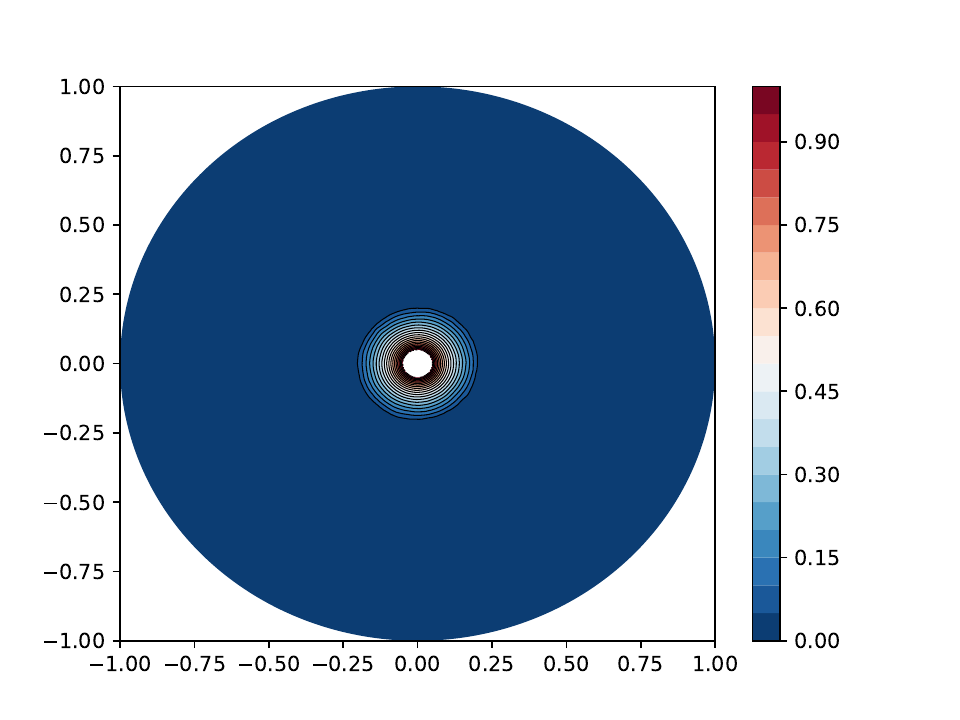}
    \end{minipage}
    \hfill%
    \begin{minipage}[c]{.46\linewidth}
        \centering
      \includegraphics[scale=1.5,width=8cm]{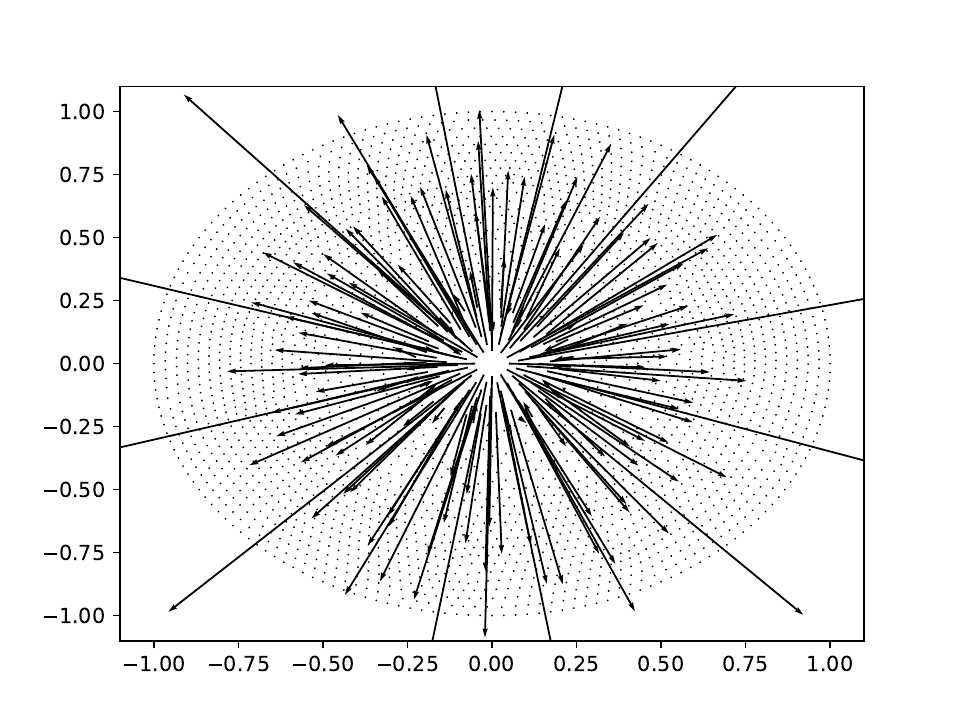}
    \end{minipage}
     \caption{An annulus  with a piececewise constant $\rho$ (discontinuous far from the interior boundary): 20 level lines (left) and a gradient map (right)\label{annulus_discontinu1}}
\end{figure}

 The second choice consists in taking a non radially symmetric coefficient $\rho$, namely
 \[
 \rho(x)=\left\{
 \begin{array}{ll}
 \rho_- & \text{ if } y<y_1,
 \\
 \rho_+& \text{ if } y> y_1,
 \end{array}
 \right.
\]
where $y_1$ is a real number in $(-1,1)$. In the case $y_1=0$, we readily check that the exact solution $u$ of \eqref{EDPinu}-\eqref{bcu} is radial and is given   in polar coordinates $(r,\theta)$ by
\[
u(x)=\frac{\ln r}{\ln 0.05}.
\]
In other words, it is independent of $\rho_-$ and $\rho_+$ and has no critical points.
If $y_1$ is different from 0, there is no analytic expression of the solution. We therefore perform a numerical simulation with  $y_1=0.5$ that corresponds to the case when the line of discontinuity of $\rho$ does not cut the interior boundary, $\rho_-=1$, and   $\rho_+=1001$. The level lines and the gradient map (again on a coarse mesh) are presented in Figure \ref{annulus_discontinu1}, and indicate that again $u$ has no critical points. Further even if the solution is singular at the interecting points between the exterior boundary and the line of discontinuity of $\rho$ (see Remark \ref{rk:reg}), the singular behavior seems to be very weak
and therefore  the influence of the line of discontinuity on the level lines is weak.

\begin{figure}[h]
\begin{minipage}[c]{.46\linewidth}
        \centering
       \includegraphics[scale=1.5,width=8cm]{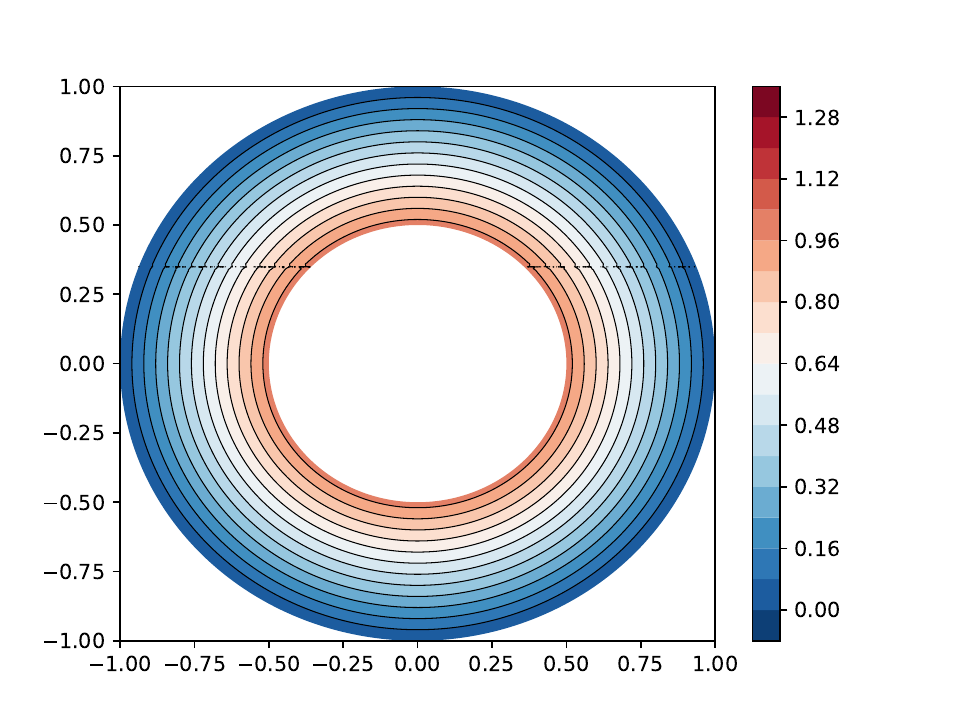}
    \end{minipage}
    \hfill%
    \begin{minipage}[c]{.46\linewidth}
        \centering
     \includegraphics[scale=1.5,width=8cm]{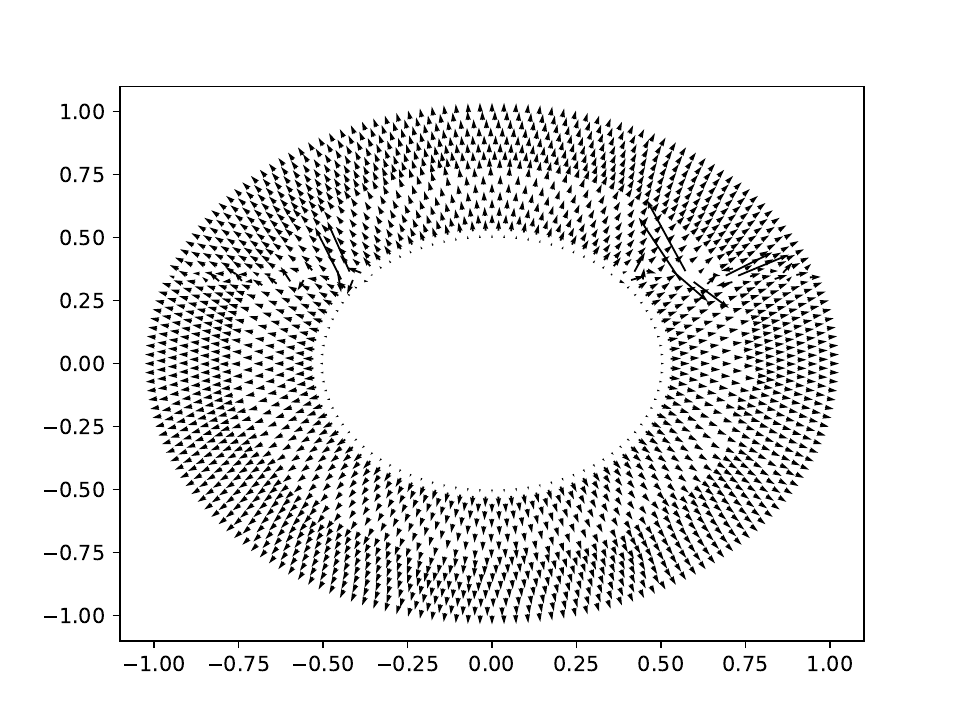}
    \end{minipage}
     \caption{An annulus  with a piececewise constant $\rho$ (discontinuous line cuts the interior boundary): 20 level lines (left) and a gradient map (right)\label{annulus_discontinu2}}
\end{figure}

Finally we consider the case where the line of discontinuity cuts the interior boundary. In that case the solution is also singular at the interecting points between the interior boundary and the line of discontinuity (see Remark \ref{rk:reg}).
For the annulus from subsection 
 \ref{ss:annulus} , by taking $y_1=0.025$ we obtain unsatisfactory numerically experiments probably due to the singular behavior of the solution and the fact that the radius of the interior ball is too small. Hence we perform a test with an interior ball of radius $0.5$ and the line of discontinuity of $\rho$ equal to $y=0.35$, $\rho_-=1$, and   $\rho_+=1001$.  From the obtained level lines and the gradient map (again on a coarse mesh) shown in Figure \ref{annulus_discontinu2}, the influence of the line discontinuity on the level lines is  much more apparent, but again it seems that $u$ has no critical points.

\section*{Acknowledgements}
R. Magnanini is partially supported by the PRIN grant n. 20229M52AS of the Italian Ministry of University (MUR) by the Gruppo Nazionale per l'Analisi Matematica, la Probabilit\`a e Applicazioni (GNAMPA) of the Istituto Nazionale di Alta Matematica (INdAM).
M. Chauvier and S. Nicaise are grateful to Prof. Christophe Troestler (Universit\'e de Mons,
  Belgium) and J. Venel (Universit\'e Polytechnique Hauts-de-France, France) for valuable discussions
on the topic of this paper.

 \appendix
 
 \section{Conformal mappings preserve  elliptic equations in divergence and diagonal form \label{sappendixA}}
 
 Let $U$ be an open and non empty set of $\R^2$
 and let $\rho\in L^\infty(U)$ be such that 
 \be\label{bdrhoappendix}
 \rho\geq 0 \hbox{ a.e. in } U.
\ee

Let $u\in H^1(U)$ be a weak solution of
\[
\div(\rho \nabla u)= 0\hbox{ in } U,
\]
in the sense that
\be\label{EDPinuappendix}
 \int_U \rho\, \nabla u\cdot \nabla v\,dx dy =0, \forall v\in H^1_0(U).
 \ee

{Let $F: U\to\mathbb{C}$ be a conformal mapping and define functions $\hat u, \hat\rho: F(U)\to\mathbb{C}$ by setting $u(z)=\hat u(F(z))$ and $\rho(z)=\hat\rho(F(z))$.  The next lemma shows that $\hat u$ is a weak solution of 
\be\label{EDPinuappendixhatPDE}
\div(\hat \rho \nabla \hat u)= 0\hbox{ in } F(U).
\ee
 \bl\label{l:preservation}
  Under the previous assumptions, $\hat u\in H^1(\hat U)$ is a weak solution of
\eqref{EDPinuappendixhatPDE},
in the sense that
\be\label{EDPinuappendixhat}
 \int_{\hat U}  \hat \rho   \nabla \hat u\cdot \hat \nabla \hat v=0, \forall \hat v\in H^1_0(\hat U).
 \ee
 \el
\begin{proof}
Let us first take $u, v\in C^1_0(U)$.
Notice that in complex notation $\nabla u\cdot \nabla v=4\,\re(\partial_{z} u\,\partial_{\bar z} v)$. By the chain rule, we have that $\partial_{z} u=(\partial_{\zeta} \hat u) \,F'(z)$, since $F$ is holomorphic. Here, we use the change of variable: $\zeta=\xi+i\,\eta=F(x+i y)=F(z)$.  Thus, we deduce that
$$
\nabla u\cdot \nabla v=4\,\re(\partial_{\zeta} \hat u\,\partial_{\bar \zeta} \hat v)\,|F'|^2=(\nabla\hat  u\cdot\nabla\hat v)\, |F'|^2.
$$
Therefore, being as $|F'|^2\,dx dy=d\xi d\eta$, we have that
$$
0=\int_U \rho\, \nabla u\cdot \nabla v\,dx dy=
\int_{F(U)} \hat\rho\, \nabla \hat u\cdot \nabla \hat v\,d\xi d\eta.
$$
 Notice that $u, v\in H_0^1(U)$ if and only if $\hat u, \hat v\in H_0^1(F(U))$. Thus, the desired conclusion follows by density from the last identity.
\end{proof}
}

        \protect\bibliographystyle{abbrv}
    \bibliography{SN-RM-MC}

@article {Seftel:63,
    AUTHOR = {\v{S}eftel , Z. G.},
     TITLE = {Estimates in {$L_{p}$} of solutions of elliptic equations
              with discontinuous coefficients and satisfying general
              boundary conditions and conjugacy conditions},
   JOURNAL = {Dokl. Akad. Nauk SSSR},
  FJOURNAL = {Doklady Akademii Nauk SSSR},
    VOLUME = {149},
    NUMBER = {},
    year ={1963},
     PAGES = {48--51},
      ISSN = {0002-3264},
   MRCLASS = {35.45},
  MRNUMBER = {196267},
MRREVIEWER = {D. E. Edmunds},
}

@article {Magnanini:80,
    AUTHOR = {Magnanini, Rolando},
     TITLE = {Sulla lipschitzianit\`a delle soluzioni di equazioni ellittiche con coefficienti radiali. ({T}he {L}ipschitz property of the solutions of elliptic
              equations with radial coefficients)},
   JOURNAL = {Le Matematiche (Catania)},
  FJOURNAL = {Le Matematiche},
    VOLUME = {35},
      YEAR = {1980},
    NUMBER = {1-2},
     PAGES = {311--324 (1983)},
      ISSN = {0373-3505},
   MRCLASS = {35J15 (35D10)},
  MRNUMBER = {698753},
MRREVIEWER = {O. John},
NOTE={(In Italian)},
}

@article(nicaise:94a,
    AUTHOR  = {S. Nicaise and A.-M. S\"andig},
    TITLE   = {General interface problems {I}},
    journal = {Math. Methods in the Appl. Sc.},
    YEAR    = 1994,
    volume  =17,
    pages={395--429}
    )

@article {MR3500302,
    AUTHOR = {Cruz-Uribe, David and Moen, Kabe and Rodney, Scott},
     TITLE = {Regularity results for weak solutions of elliptic {PDE}s below
              the natural exponent},
   JOURNAL = {Ann. Mat. Pura Appl. (4)},
  FJOURNAL = {Annali di Matematica Pura ed Applicata. Series IV},
    VOLUME = {195},
      YEAR = {2016},
    NUMBER = {3},
     PAGES = {725--740},
      ISSN = {0373-3114},
   MRCLASS = {35J25 (35B45 35B65 35D30 42B37 46E35)},
  MRNUMBER = {3500302},
       DOI = {10.1007/s10231-015-0486-y},
       URL = {https://doi-org.proxybmath.univ-lille.fr/10.1007/s10231-015-0486-y},
}

@article {Meyers:63,
    AUTHOR = {Meyers, Norman G.},
     TITLE = {An {$L^{p}$}e-estimate for the gradient of solutions of
              second order elliptic divergence equations},
   JOURNAL = {Ann. Scuola Norm. Sup. Pisa Cl. Sci. (3)},
  FJOURNAL = {Annali della Scuola Normale Superiore di Pisa. Classe di
              Scienze. Serie III},
    VOLUME = {17},
      YEAR = {1963},
     PAGES = {189--206},
      ISSN = {0391-173X},
   MRCLASS = {35.42},
  MRNUMBER = {159110},
MRREVIEWER = {Lipman Bers},
}

@article{Perelmuter:25,
author = {Perelmuter, M.},
year = {2025},
month = {04},
pages = {},
title = {Second order regularity of solutions of elliptic equations in divergence form with Sobolev coefficients},
journal = {Annali di Matematica Pura ed Applicata (1923 -)},
doi = {10.1007/s10231-025-01569-w}
}

@article {PiccininiSpagnolo:72,
    AUTHOR = {Piccinini, L. C. and Spagnolo, S.},
     TITLE = {On the {H}\"{o}lder continuity of solutions of second order
              elliptic equations in two variables},
   JOURNAL = {Ann. Scuola Norm. Sup. Pisa Cl. Sci. (3)},
  FJOURNAL = {Annali della Scuola Normale Superiore di Pisa. Classe di
              Scienze. Serie III},
    VOLUME = {26},
      YEAR = {1972},
     PAGES = {391--402},
      ISSN = {0391-173X},
   MRCLASS = {35J15},
  MRNUMBER = {361422},
MRREVIEWER = {C. A. Swanson},
}

@Misc{Netgen,
  author =   {Sch\"{o}berl, Joachim},
  title =    {Netgen/NGSolve},
  howpublished = {\texttt{https://ngsolve.org/}},
  year =     {2022}}

@article {FabesKenigSerapioni:82,
    AUTHOR = {Fabes, Eugene B. and Kenig, Carlos E. and Serapioni, Raul P.},
     TITLE = {The local regularity of solutions of degenerate elliptic
              equations},
   JOURNAL = {Comm. Partial Differential Equations},
  FJOURNAL = {Communications in Partial Differential Equations},
    VOLUME = {7},
      YEAR = {1982},
    NUMBER = {1},
     PAGES = {77--116},
      ISSN = {0360-5302},
   MRCLASS = {35J70 (35D10)},
  MRNUMBER = {643158},
MRREVIEWER = {M.-T. Lacroix},
       DOI = {10.1080/03605308208820218},
       URL = {https://doi-org.proxybmath.univ-lille.fr/10.1080/03605308208820218},
}

@article {DeCiccoVivaldi:96,
    AUTHOR = {De Cicco, Virginia and Vivaldi, Maria Agostina},
     TITLE = {Harnack inequalities for {F}uchsian type weighted elliptic
              equations},
   JOURNAL = {Comm. Partial Differential Equations},
  FJOURNAL = {Communications in Partial Differential Equations},
    VOLUME = {21},
      YEAR = {1996},
    NUMBER = {9-10},
     PAGES = {1321--1347},
      ISSN = {0360-5302},
   MRCLASS = {35J70 (35B45 35J15)},
  MRNUMBER = {1410832},
MRREVIEWER = {Gaston E. Hernandez},
       DOI = {10.1080/03605309608821229},
       URL = {https://doi-org.proxybmath.univ-lille.fr/10.1080/03605309608821229},
}

@article {Alessandrini:87,
    AUTHOR = {Alessandrini, Giovanni},
     TITLE = {Critical points of solutions of elliptic equations in two
              variables},
   JOURNAL = {Ann. Scuola Norm. Sup. Pisa Cl. Sci. (4)},
  FJOURNAL = {Annali della Scuola Normale Superiore di Pisa. Classe di
              Scienze. Serie IV},
    VOLUME = {14},
      YEAR = {1987},
    NUMBER = {2},
     PAGES = {229--256 (1988)},
}

@article {Alessandrini:87b,
    AUTHOR = {Alessandrini, Giovanni},
     TITLE = {An identification problem for an elliptic equation in two
              variables},
   JOURNAL = {Ann. Mat. Pura Appl. (4)},
  FJOURNAL = {Annali di Matematica Pura ed Applicata. Serie Quarta},
    VOLUME = {145},
      YEAR = {1986},
     PAGES = {265--295},
      ISSN = {0003-4622},
   MRCLASS = {35R30 (35J25)},
  MRNUMBER = {886713},
MRREVIEWER = {S. D\"{u}mmel},
       DOI = {10.1007/BF01790543},
       URL = {https://doi-org.proxybmath.univ-lille.fr/10.1007/BF01790543},
}

@article {AlessandriniLupoRosset:93,
    AUTHOR = {Alessandrini, Giovanni and Lupo, Daniela and Rosset, Edi},
     TITLE = {Local behavior and geometric properties of solutions to
              degenerate quasilinear elliptic equations in the plane},
   JOURNAL = {Appl. Anal.},
  FJOURNAL = {Applicable Analysis. An International Journal},
    VOLUME = {50},
      YEAR = {1993},
    NUMBER = {3-4},
     PAGES = {191--215},
      ISSN = {0003-6811,1563-504X},
   MRCLASS = {35J60 (35Bxx)},
  MRNUMBER = {1278325},
MRREVIEWER = {Yao\ Tian\ Shen},
       DOI = {10.1080/00036819308840193},
       URL = {https://doi.org/10.1080/00036819308840193},
}

@article {AlessandriniMagnanini:92,
    AUTHOR = {Alessandrini, G. and Magnanini, R.},
     TITLE = {The index of isolated critical points and solutions of
              elliptic equations in the plane},
   JOURNAL = {Ann. Scuola Norm. Sup. Pisa Cl. Sci. (4)},
  FJOURNAL = {Annali della Scuola Normale Superiore di Pisa. Classe di
              Scienze. Serie IV},
    VOLUME = {19},
      YEAR = {1992},
    NUMBER = {4},
     PAGES = {567--589},
      ISSN = {0391-173X},
   MRCLASS = {35J25 (35B05 35B30 58E05)},
  MRNUMBER = {1205884},
MRREVIEWER = {Xin Kang Guo},
       URL = {http://www.numdam.org.proxybmath.univ-lille.fr/item?id=ASNSP_1992_4_19_4_567_0},
}

@article {AlessandriniMagnanini:94,
    AUTHOR = {Alessandrini, G. and Magnanini, R.},
     TITLE = {Elliptic equations in divergence form, geometric critical
              points of solutions, and {S}tekloff eigenfunctions},
   JOURNAL = {SIAM J. Math. Anal.},
  FJOURNAL = {SIAM Journal on Mathematical Analysis},
    VOLUME = {25},
      YEAR = {1994},
    NUMBER = {5},
     PAGES = {1259--1268},
      ISSN = {0036-1410},
   MRCLASS = {35P15 (35J25)},
  MRNUMBER = {1289138},
MRREVIEWER = {Jan Bochenek},
       DOI = {10.1137/S0036141093249080},
       URL = {https://doi-org.proxybmath.univ-lille.fr/10.1137/S0036141093249080},
}

@unpublished{Chauvieretco:2025,
  TITLE = {An existence result and simulations of a space charges problem applied to {H}{V}{D}{C} transmission lines},
  AUTHOR = {Chauvier, Madeline and  Nicaise, Serge and  Troestler, Christophe
  and Venel, Juliette},
  NOTE = {in preparation},
  YEAR = {2025},
}

@article {Dengetco:18,
    AUTHOR = {Deng, Haiyun and Liu, Hairong and Tian, Long},
     TITLE = {Critical points of solutions to a quasilinear elliptic
              equation with nonhomogeneous {D}irichlet boundary conditions},
   JOURNAL = {J. Differential Equations},
  FJOURNAL = {Journal of Differential Equations},
    VOLUME = {265},
      YEAR = {2018},
    NUMBER = {9},
     PAGES = {4133--4157},
      ISSN = {0022-0396},
   MRCLASS = {35J60 (35B38 35J25 35J93)},
  MRNUMBER = {3843296},
       DOI = {10.1016/j.jde.2018.05.031},
       URL = {https://doi-org.proxybmath.univ-lille.fr/10.1016/j.jde.2018.05.031},
}

@article {Dengetco:22,
    AUTHOR = {Deng, Haiyun and Liu, Hairong and Yang, Xiaoping},
     TITLE = {Critical points of solutions to a kind of linear elliptic
              equations in multiply connected domains},
   JOURNAL = {Israel J. Math.},
  FJOURNAL = {Israel Journal of Mathematics},
    VOLUME = {249},
      YEAR = {2022},
    NUMBER = {2},
     PAGES = {935--971},
      ISSN = {0021-2172},
   MRCLASS = {35B38},
  MRNUMBER = {4462650},
       DOI = {10.1007/s11856-022-2330-6},
       URL = {https://doi-org.proxybmath.univ-lille.fr/10.1007/s11856-022-2330-6},
}

@article {FinnGilbarg:57,
    AUTHOR = {Finn, R. and Gilbarg, D.},
     TITLE = {Asymptotic behavior and uniquenes of plane subsonic flows},
   JOURNAL = {Comm. Pure Appl. Math.},
  FJOURNAL = {Communications on Pure and Applied Mathematics},
    VOLUME = {10},
      YEAR = {1957},
     PAGES = {23--63},
      ISSN = {0010-3640},
   MRCLASS = {76.0X},
  MRNUMBER = {86556},
MRREVIEWER = {M. J. Lighthill},
       DOI = {10.1002/cpa.3160100102},
       URL = {https://doi-org.proxybmath.univ-lille.fr/10.1002/cpa.3160100102},
}

@article {Magnanini:16,
    AUTHOR = {Magnanini, Rolando},
     TITLE = {An introduction to the study of critical points of solutions
              of elliptic and parabolic equations},
   JOURNAL = {Rend. Istit. Mat. Univ. Trieste},
  FJOURNAL = {Rendiconti dell'Istituto di Matematica dell'Universit\`a di
              Trieste. An International Journal of Mathematics},
    VOLUME = {48},
      YEAR = {2016},
     PAGES = {121--166},
      ISSN = {0049-4704},
   MRCLASS = {35B38 (35J05 35J08 35J10 35J15 35K10)},
  MRNUMBER = {3592440},
       DOI = {10.13137/2464-8728/13154},
       URL = {https://doi-org.proxybmath.univ-lille.fr/10.13137/2464-8728/13154},
}

@book(gilbarg:77,
   author = {D.~Gilbarg and N.~S.~Trudinger},
   title  = {{E}lliptic {P}artial {D}ifferential {E}quations of {S}econd
             {O}rder},
   year   = 1977,
   publisher = {Springer},
   address = {Berlin -- Heidelberg -- New York})

@article {SabinadeLis:15,
    AUTHOR = {Sabina de Lis, Jos\'{e} C.},
     TITLE = {Hopf maximum principle revisited},
   JOURNAL = {Electron. J. Differential Equations},
  FJOURNAL = {Electronic Journal of Differential Equations},
      YEAR = {2015},
     PAGES = {No. 115, 9},
   MRCLASS = {35B50 (35J25)},
  MRNUMBER = {3358487},
}

@book{AmorusoLattarulo:14,
   title =     {Filamentary ion flow : theory and experiments},
   author =    {Amoruso, Vitantonio; Lattarulo, Francesco},
   publisher = {Wiley-IEEE Press},
  year = {2014},
month = {01},
pages = {1-204},}

@article{GuillodPfeifferFranck:14,
author = {Guillod, Thomas and Pfeiffer, Michael and Franck, Christian},
year = {2014},
month = {12},
pages = {2493-2501},
title = {Improved Coupled Ion-Flow Field Calculation Method for AC/DC Hybrid Overhead Power Lines},
volume = {29},
journal = {Power Delivery, IEEE Transactions on},
doi = {10.1109/TPWRD.2014.2322052}
}

@article{Xiaoetall:17,
author = {Xiao, Fengnyu and Zhang, Bo and Mo, Jianghua and Jinliang, He},
year = {2017},
month = {08},
pages = {1-1},
title = {Calculation of 3-D Ion Flow Field at the Crossing of HVDC Transmission Lines by Method of Characteristics},
volume = {PP},
journal = {IEEE Transactions on Power Delivery},
doi = {10.1109/TPWRD.2017.2737006}
}

@article{Zhangetall:15,
author = {Zhang, Bo and Mo, Jianghua and Yin, Han and Jinliang, He},
year = {2015},
month = {03},
pages = {1-4},
title = {Calculation of Ion Flow Field Around HVdc Bipolar Transmission Lines by Method of Characteristics},
volume = {51},
journal = {Magnetics, IEEE Transactions on},
doi = {10.1109/TMAG.2014.2349491}
}

@article{ButlerCendesHoburg:89,
author = {Butler, A.J. and Cendes, Z.J. and Hoburg, J.F.},
year = {1989},
month = {06},
pages = {533 - 538},
title = {Interfacing the finite-element method with the method of characteristics in self-consistent electrostatic field models},
volume = {25},
journal = {Industry Applications, IEEE Transactions on},
doi = {10.1109/28.31225}
}

@article{CristinaDinelliFeliziani:91,
author = {Cristina, Sebastian and Dinelli, Giorgio and Feliziani, Mauro},
year = {1991},
month = {02},
pages = {147 - 153},
title = {Numerical computation of corona space charge and V-I characteristic in DC electrostatic precipitators},
volume = {27},
journal = {Industry Applications, IEEE Transactions on},
doi = {10.1109/28.67546}
}

@article{Lobry:14,
author = {Lobry, Jacques},
year = {2014},
month = {02},
pages = {541-544},
title = {A New Numerical Scheme for the Simulation of Corona Fields},
volume = {50},
journal = {Magnetics, IEEE Transactions on},
doi = {10.1109/TMAG.2013.2281999}
}

@article{SharmaJaniszewski:69,
author = {Sharma, M. P. and  Janiszewski, V. },
year = {1969},
month = {05},
pages = {718--731},
title = {Analysis of Corona Losses on DC transmission Lines: I-Unipolar Lines},
volume = {88},
journal = {IEEE Trans. Power App. Syst.},
}

\end{document}